\documentclass[11pt]{amsart}
\usepackage{amssymb}
\usepackage{times}
\usepackage{url}
\usepackage{hyperref}

\hypersetup{pdfpagemode=UseNone}

\newtheorem{theorem}{Theorem}[section]
\newtheorem{lemma}[theorem]{Lemma}
\newtheorem{proposition}[theorem]{Proposition}
\newtheorem{corollary}[theorem]{Corollary}
\theoremstyle{remark}
\newtheorem{remark}[theorem]{Remark}

\theoremstyle{example}
\newtheorem{example}[theorem]{Example}

\newcommand{\F}{\mathbb{F}}

\newcommand{\eL}{\mathbb{L}}
\newcommand{\E}{\mathbb{E}}
\newcommand{\Z}{\mathbb{Z}}
\newcommand{\Q}{\mathbb{Q}}
\newcommand{\C}{\mathbb{C}}
\newcommand{\R}{\mathbb R}
\renewcommand{\P}{\mathbb{P}}
\newcommand{\GL}{\mathrm{GL}}
\newcommand{\SL}{\mathrm{SL}}
\newcommand{\Sp}{\mathrm{Sp}}

\newcommand{\rk}{\mathrm{rk}}
\newcommand{\hi}{\mathrm{h}}
\newcommand{\UT}{\mathrm{UT}}

\newcommand{\glnz}{\mathrm{GL}(n,\Z)}

\newcommand{\glnf}{\mathrm{GL}(n,\F)}

\newcommand{\cha}{\mathrm{char} \hspace{.5pt}}

\newcommand{\calg}{{\mathcal G}}

\def\gp#1{\langle \hspace*{.2mm} #1 \hspace*{.25mm} \rangle}

\newcommand{\gpess}{\langle S \hspace{1pt} \rangle}
\newcommand{\st}{\hspace{2.5pt} | \hspace{2.5pt}}
\newcommand{\abk}{\allowbreak}

\begin{document}

\title{Linear groups and computation}

\begin{abstract}
We present an exposition of our ongoing project in a new
area of applicable mathematics: practical computation with 
finitely generated linear groups over infinite fields. 
Methodology and algorithms available for this class of groups 
are surveyed. We 
illustrate the solution of hard 
problems by
computer experimentation. Possible avenues for further progress 
are discussed. 
\end{abstract}

\author{A.~S.~Detinko}

\author{D.~L.~Flannery}



\maketitle

\section{Introduction}

A \emph{linear group} (interchangeably, \emph{matrix group}) of 
degree $n$ over a field $\F$ is a subgroup of the group $\glnf$ 
of all $n\times n$ invertible matrices with entries in $\F$. 
Linear groups have impact throughout mathematics, in diverse 
branches such as topology, geometry, number theory, and analysis. 
They are amenable to calculation and modeling 
transformations, so provide a bridge between algebra and the
physical sciences. 
For example, applications arise in crystallography~\cite{AClib}, 
theoretical physics~\cite{MonodromyAppendix}, 
error-correcting codes~\cite{Zinoviev}, 
cryptography~\cite{ShpilrainBook},  
even quantum computing~\cite{Derksen}.

The study of matrix groups goes back to the origin of group theory.
Klein and Lie discovered their role in geometry and 
differential equations, but the foundations were really laid
down by Jordan~\cite{Jordan}.
Pioneering contributions were also made by Schur, 
Burnside, Frobenius, and Blichfeldt.

Linear group theory supplements representation theory. 
It has played a large part in the classification of finite 
simple groups (see \cite{Aschbacher84},
\cite[pp.~76--78]{FeitSurvey}) and the theory of infinite 
solvable groups~\cite[Chapter 3]{LennoxRobinson}. 
Important classes of groups exhibit linearity: e.g., 
each polycyclic-by-finite or countable free group
has a faithful representation over the 
integers $\Z$. 

Detailed accounts of our subject may be found in the books 
\cite{DixonBook,SuprunenkoI,WBook} 
and surveys \cite{ZalesskiiI,ZalesskiiIII,ZalesskiiII}.

\subsection{CGT and matrix groups}

Computational group theory (CGT) is a modern discipline interfacing 
algebra and computer science. It deals with 
the design, implementation, and analysis of algorithms for 
groups and related objects.  Progress in CGT is tied to the 
evolution of computer algebra systems such as 
{\sf GAP}~\cite{GAP}, {\sc Magma}~\cite{Magma}, and 
SageMath~\cite{Sage}. 
Other systems (e.g., Maple~\cite{Maple} and 
Mathematica~\cite{Mathematica}) allow some group-theoretical 
computation. We refer to 
\cite{HoltBook,HulpkeNotes,AkosSurvey,Sims} for 
coverage of aspects of CGT.

Typically a group is input  to the computer as permutations,  
or as a presentation (generating elements and relations 
between them), or as matrices acting on a vector space. Matrix 
representations have the advantage of compactness; a huge 
(even infinite) group can be represented by input of small 
size. For example, the smallest permutation representation 
of the Monster simple group has degree 
about $10^{20}$, whereas R.~A.~Wilson found a linear representation 
of the Monster in dimension $196882$ over the $2$-element field. 
(An early triumph of CGT was its deployment in 
proofs of existence of sporadic simple groups~\cite[p.~4]{HoltBook}.)
At the other end of the historical spectrum, Jordan's
description of solvable matrix groups over a finite field 
can be viewed as a prototype CGT algorithm, constructing groups 
inductively from ones in smaller degrees.

Matrix groups have become a primary focus of CGT.
Currently the emphasis is on practical algorithms that
permit fruitful computer experimentation.
Tremendous effort has been expended on groups over finite fields, 
coalescing around the `Matrix Group Recognition 
Project'~\cite{BBS,OBriensurvey}. This is now at an advanced stage, 
thanks to many authors. It comprises a substantial part of the armory 
of CGT.  Given a matrix group $G$ over a finite field,
preliminary objectives are to compute the composition series of $G$ and 
recognize series factors. Thereafter one determines a 
presentation of $G$, and proceeds to tasks such as exploring
subgroup structure.

\subsection{Computing with linear groups over infinite domains}

Matrix groups over infinite 
fields have received 
comparatively scant attention in CGT. Since these groups are potentially
infinite, methods for groups over finite fields or for permutation groups 
are not immediately applicable. Note that infinite linear groups need 
not be finitely presentable. Even if we have a presentation, we may still 
be unable to use it or other helpful knowledge about the 
group efficiently. Complexity issues (such as growth of matrix 
entries during computation) cause bottlenecks.
More seriously, certain algorithmic problems for infinite groups
are undecidable: an algorithm to solve the problem for all input does 
not exist (see \cite[Section~5]{MillerSurvey}). 

Despite the pitfalls noted above, there is reason to be optimistic. 
In \cite{BBR,Beals1,Beals2}, Babai and his collaborators initiated the 
development of algorithms for linear groups over infinite fields 
(mainly the rationals $\Q$). They rely on algorithms for matrix 
algebras; cf.~\cite{Ronyai}.  
These papers 
establish 
existence of low complexity (e.g., polynomial-time) solutions of
the algorithmic problems considered.
Several years later, people became interested in 
computing with linear groups of special kinds.  Representations 
over $\Z$ figure prominently in the 
theory of polycyclic 
groups (\cite{Baumslagetal}, \cite[Section 9.2]{LennoxRobinson}). 
This is a gateway to algorithms for polycyclic groups where the
input is a finite set of generating 
matrices over $\Z$ or $\Q$. Inspired by Dixon~\cite{Dixon85}, work 
in this area was undertaken by Lo and 
Ostheimer~\cite{OstheimerLo,Ostheimer}, then Assmann and 
Eick~\cite{BettinaBjorn,BettinaBjorn2}.
The latter was incorporated into \cite{Polenta}, one of the 
first software packages dedicated to computing with infinite 
linear groups.

On a separate front, algorithms for linear algebraic groups have been 
developed by Cohen et al. beginning in the late 1990s~\cite{Cohenetal,gra13}.  
Here the input groups are not finitely generated, so are given by a finite 
set of polynomials, or Lie algebras, rather than a finite set of matrices.
See \cite{deGraafPrograms,deGraafBook,WillemBook} for  
pertinent material here.
 
A massive amount of classification data for linear groups over infinite
fields has accumulated 
over the years. While much of this is yet to be implemented, some notable 
classifications are available in {\sf GAP} and {\sc Magma}. One 
of these is the library of almost crystallographic groups in 
\cite{AClib}. 
Another is the {\sc Magma} database of irreducible maximal finite subgroups of 
$\glnz$ for $n\leq 23$ derived from work by Plesken et al. on 
integral representation of finite groups (see also~\cite{NebePlesken}). 

Several cognate topics lie outside the scope of this article. 
Perhaps the most relevant is computational representation theory, 
which has been active from the 1960s up to now.

\subsection{Algorithms for finitely generated linear groups}

To compute with linear groups over infinite fields, 
we must first decide how the groups will be designated 
in a computer. Of course, input will be a finite set. A
convenient format is a finite set of generating matrices. 
Although not every matrix group can be so designated, 
those that are constitute a major 
class,
and occur frequently in applications. 

This article is an exposition of our ongoing project to
compute with groups given by a finite set of generating matrices
over an (arbitrary) infinite field. We

\begin{itemize}
\item[(i)] formulate general methodology;
\item[(ii)] apply (i) to design effective algorithms;
\item[(iii)] implement the algorithms and demonstrate their 
practicality.
\end{itemize}
As (ii) and (iii) indicate, an overarching goal is to obtain 
algorithms that complete in reasonable time for a wide range of 
inputs. Ideally, the software would replace traditional 
mathematics by machine computation, simplifying the 
solution of problems, and leading to the solution of formerly 
intractable problems.

Our methodology draws on (classical) theory of 
linear groups, as in \cite{DixonBook,WBook}. This equips us 
with well-tried tools, 
such as the `method of finite approximation'. 
Apart from underpinning the success of our approach, 
linear group theory and its central concerns 
guide our choice of problems to 
give priority. One of these is realizing the 
Tits alternative computationally. That is, we 
devise and implement 
a practical algorithm to test whether a 
finitely generated linear group is solvable-by-finite 
(recall that $G$ is \emph{${\sf X}$-by-finite}, 
or \emph{virtually ${\sf X}$}, if there exists a finite index normal 
subgroup of $G$ that has property ${\sf X}$).
Then we dispose of further questions for solvable-by-finite linear groups: 
recognition problems such as 
testing whether a group is finite, solvable, or nilpotent. 
Later parts of the article are occupied with the second class of 
the Tits alternative, specifically arithmetic 
and Zariski dense subgroups of semisimple algebraic 
groups. We conclude by discussing avenues for future research.

The aims of this project were achieved with
our colleagues Willem de Graaf, Bettina Eick, Alexander Hulpke, 
and Eamonn O'Brien, to whom we are deeply grateful.

\medskip

\section{Computing with finitely generated linear groups}

This section introduces some of the basic ideas in 
computing with  linear groups by means of 
congruence homomorphisms. 

\subsection{How to input groups}\label{RepGps}

The computer algebra system in which we implement our
algorithms must support computation 
over the defining field or ring of the input.
A field that one often meets in practice is $\R$.  
Floating point representation of real numbers is popular in 
applied mathematics, but it is unsuitable for computing with 
linear groups defined over an arbitrary field $\F$: CGT 
entails inherently symbolic computation; output should be exact. 

Let $G = \langle S \rangle$ where 
$S= \{ g_1,\ldots ,g_r\}\subseteq \glnf$.
Suppose for simplicity that $g^{-1}\in S$ if $g\in S$.
The group $G$ is defined over the subfield of $\F$ 
generated by all entries of the $g_i$. So we assume that
$\F$ is a finitely generated extension of its prime subfield $\F_0$. 
\begin{lemma}\label{Lemma1}
There exist algebraically independent ${\rm x}_1, \ldots, {\rm x}_m\in \F$, 
$m\geq 0$, such that $\F$ is a finite extension of the function field
$\F_0({\rm x}_1, \ldots, {\rm x}_m)$.
\end{lemma}
\begin{remark} For an algorithm to compute  
the ${\rm x}_i$ as in Lemma~\ref{Lemma1}, see \cite{Steinwandt}.
\end{remark}

By Lemma~\ref{Lemma1}, `arbitrary field' for us is one of 

\begin{itemize}
\item[(I)]  $\Q$;
\item[(II)] a number field $\P$ (finite degree extension of $\Q$);
\item[(III)] a function field $\E({\rm x}_1,\ldots , {\rm x}_m)$, $m\geq 1$, 
where $\E = \Q$, $\P$, or
a finite field $\F_q$ of size $q$;
\item[(IV)] a finite degree extension of $\E({\rm x}_1,\ldots , {\rm x}_m)$.
\end{itemize}

\medskip

\noindent Happily, {\sc Magma} supports computation in these fields.
Of course, (I) and (III) are specializations of (II) and (IV), respectively; 
but we choose to distinguish between the four 
types. When computing over number fields 
$\P$, we do not inflate matrix dimension by the degree $|\P:\Q|$ according
to the action of $\P$ on a $\Q$-basis of $\P$, to
get matrices with entries in $\Q$. This is a standard trick, but
risks increasing the dimension beyond the boundaries of practicality.

\subsection{Finite approximation}\label{FASubsection}
 
Let $R$ be the subring of $\F$ generated by the entries of all elements 
of $S$. Hence $R$ is a finitely generated integral domain and $G\leq \abk
\GL(n,R)$. 
\begin{lemma}[{\cite[p.~50]{WBook}}] 
\label{WFGID}
The ring $R$ is approximated by finite fields; i.e., 
\begin{itemize}
\item[{\rm (i)}] if $\varrho$ is a maximal ideal of $R$
 then $R/\varrho$ is a finite field,
\item[{\rm (ii)}] if $a\in R\setminus \{ 0\}$ then $R$ has a maximal ideal 
not containing
$a$.
\end{itemize}
\end{lemma}

Lemma~\ref{WFGID} underlies the following, 
due to Mal'cev (see \cite[(2.1), p.~74]{ZalesskiiI}).
\begin{theorem}
\label{Malcev}
Each finitely generated subgroup of $\glnf$ is 
approximated by matrix groups of degree $n$ 
over finite fields. 
\end{theorem}

Theorem~\ref{Malcev} is at the heart of our computational strategy. 
The word `approximated' expresses the residual finiteness of $G$: 
if $g\in G\setminus\{ 1\}$ then there exists 
a homomorphism $f$ from 
$G$ onto a subgroup of $\GL(n,q)$ for some prime power $q$, 
such that $f(g)\neq 1$. In the proof of Theorem~\ref{Malcev} in
\cite[Theorem~4.2, p.~51]{WBook},
$f$ is a \emph{congruence homomorphism}.
We now begin to set up the formalism that we adopt to compute with
these homomorphisms.

Let $\varrho$ be a (proper) ideal of an associative unital ring $\Delta$.
Natural surjection 
$\Delta \rightarrow \Delta/\varrho$ induces an algebra homomorphism 
$\mathrm{Mat}(n,\Delta) \rightarrow \mathrm{Mat}(n,\Delta/\varrho)$, which
then restricts to a group homomorphism 
$\mathrm{GL}(n,\Delta) \rightarrow \mathrm{GL}(n,\Delta/\varrho)$. We 
denote all these homomorphisms by $\varphi_\varrho$. The 
\emph{principal congruence subgroup (PCS) $\Gamma_{n,\varrho}$ of level $\varrho$} 
is the kernel of $\varphi_\varrho$ in $\GL(n,\Delta)$. Theorem~\ref{Malcev} tells
us that for each $g\neq 1$ in $G$ there exists a maximal ideal 
$\varrho$ of $R$ such that $\varphi_{\varrho}(g)$ is a non-identity element of 
the general linear group $\GL(n,R/\varrho)$ over the finite field $R/\varrho$.
The accuracy of information about  $G$ provided by its
congruence images $\varphi_\varrho(G)$, and how we compute 
with matrix groups over finite rings, govern the effectiveness of 
this overall approach to computing with finitely generated subgroups of $\glnf$.

\subsection{Constructing congruence homomorphisms}\label{ConstrCH}

We explain how to construct congruence homomorphisms for the  
main field types (I)--(IV), in line with Theorem~\ref{Malcev}. 

\subsubsection{}\label{VarphiForRationals}
If $\F = \Q$ then $R =  
\frac{1}{\mu}\Z= \Z[\frac{1}{\mu}]$, the ring of rationals whose denominators are 
powers of a fixed integer $\mu$. We can take $\mu$ to be the least common multiple
of the denominators of the entries of the $g_i$. For any positive 
$b\in \Z$ not dividing $\mu$, 
entrywise reduction of the $g_i$ modulo $b$ defines a congruence homomorphism 
$\varphi_\varrho$ where $\varrho=bR$. We write 
`$b$' in place of  
`$\varrho$' as a subscript in the notation.
If $b=p$ is prime then $\varrho$ is a maximal ideal of 
$R$ and $R/\varrho=\F_p$. Set $\varphi_{1,p}:=\varphi_p$.

\subsubsection{} Let $\F$ be a number field  $\P$ of degree $k$ over $\Q$.
There is an algebraic integer $\alpha$ with minimal polynomial 
$f(t) = a_0+a_1t+\cdots + a_{k-1}t^{k-1}+t^k\in \Z[t]$ such that 
$\P = \Q(\alpha)$. We have $R \subseteq \frac{1}{\mu}\Z[\alpha]$ for 
some $\mu\in \Z$. Let $p\in \Z$ be a prime not dividing $\mu$, and set 
$\bar{b}_i = \abk \varphi_p(b_i)$, 
$\bar{a}_i = \varphi_p(a_i)$. For any root $\bar{\alpha}$ of 
$\bar{f}(t)=\abk \sum_{i=0}^{k-1}\bar{a}_it^i$,
define the homomorphism $\varphi_{2,p}$ of $R$ onto the finite 
field $\F_p(\bar{\alpha})$ by
\[
\textstyle{\sum_{i=0}^{k-1} b_i\alpha^i \mapsto 
\sum_{i=0}^{k-1} \bar{b}_i\bar{\alpha}^i}.
\]
Let $\bar{f}_j(t)$ be an irreducible factor of $\bar{f}(t)$. If
 $f_j(t)$ is a preimage of $\bar{f}_j(t)$  in $\Z[t]$ then
the ideal of $R$ generated by $p$ and $f_j(\alpha)$ is maximal, and 
$\varphi_{2,p} = \varphi_{\varrho}$.

\subsubsection{}\label{FirstFunctionFieldsSubsection}
Let $\F = \E({\rm x}_1, \ldots, {\rm x}_m)$. Then
$R \subseteq \frac{1}{\mu}\E[{\rm x}_1, \ldots, {\rm x}_m]$ for some 
polynomial $\mu = \mu({\rm x}_1, \ldots, {\rm x}_m) \in 
\E [{\rm x}_1, \ldots, {\rm x}_m ]$. 
Let $\alpha = (\alpha_1,\ldots , \alpha_m)$ be a non-root of $\mu$, where 
each $\alpha_i$ is in the algebraic closure $\overline{\E}$ of $\E$. The 
$\alpha_i$ can be chosen in $\E$ if 
$\mathrm{char} \, \E = 0$ and 
in a finite extension $\F_{q^c}$ if $\E = \F_q$.
Define $\varphi_{3,\alpha}$ on $R$ 
to be the homomorphism that substitutes 
$\alpha_i$ for ${\rm x}_i$, $1\leq i\leq m$. 
Then 
$\varphi_{3,\alpha,p}:= \varphi_{i,p} \circ \varphi_{3,\alpha}$
where $i=1$ if $\E = \Q$ and $i=2$ if $\E\neq \Q$ is a number field.

\subsubsection{}\label{VarphiForAFF}

Suppose that $\F$ is an extension of $\eL = \E({\rm x}_1, \ldots, {\rm x}_m)$ 
of degree $e$: say $\F = \abk \eL(\beta)$. Such $\beta$ exist if
$\cha \, \F = 0$  or $\F$ is perfect in positive characteristic.
Let $f(t) = t^e+ \abk a_{e-1}t^{e-1}+ \cdots +a_1t+a_0$ be the minimal 
polynomial of $\beta$. Then $R\subseteq \frac{1}{\mu}\eL_0[\beta]$ where 
$\mu \in \eL_0 =  \E[{\rm x}_1, \ldots, {\rm x}_m]$. 
We may assume that $f(t)\in \eL_0[t]$.

Define $\varphi_{4,\alpha}$ on $\GL(n,R)$ as follows. 
Take a non-root $\alpha = (\alpha_1,\dots , \alpha_m)\in \overline{\E}^m$
of $\mu$. 
(Remember that we can find $\alpha_i\in \E$ if $\cha \, \E = 0$ and 
$\alpha_i \in \F_{q^c}$ if $\E = \F_q$.)
Further, let $\tilde{\beta}$ be a root of 
$\tilde{f}(t) = \sum_{i=0}^{e-1}\tilde{a}_it^i + t^e$, where 
$\tilde{a}_i=\varphi_{3,\alpha}(a_i)$.
Each element of $\frac{1}{\mu}\eL_0[\beta]$ has a unique expression as 
$\sum_{i=0}^{e-1}c_i\beta^i$ for $c_i \in \frac{1}{\mu}\eL_0$. Then
\[
\varphi_{4,\alpha}: {\textstyle \sum_{i=0}^{e-1}c_i\beta^i\mapsto 
\sum_{i=0}^{e-1}\tilde{c}_i\tilde{\beta}^i}
\]
where $\tilde{c}_i = \varphi_{3,\alpha}(c_i)$. In zero characteristic, 
$\varphi_{4,\alpha, p} := \varphi_{i,p}\circ\varphi_{4,\alpha}$ 
with $i=1$ when $\E=\Q$ and $\tilde{\beta}\in \Q$, and $i=2$ when $\E=\P$ .

\bigskip

As \ref{VarphiForRationals}--\ref{VarphiForAFF} show, constructing  
congruence homomorphisms is straightforward. 
The main operations are reduction modulo primes $p\in \Z$ and substitution 
of indeterminates. These define an appropriate 
maximal ideal $\varrho$ and hence image field.
The sort of ideal $\varrho$ chosen will depend on the problem at hand.

\subsection{Computational finite approximation}

We continue with the notation above: $R\subseteq \F$ is determined 
by a generating set for $G\leq \glnf$, $\varrho$ is a maximal ideal of $R$, 
and $\varphi_\varrho$ is the corresponding homomorphism 
$\GL(n,R)\rightarrow \GL(n,R/\varrho)$, 
with $R/\varrho$ a finite field. Denote the kernel of $\varphi_\varrho$ 
on $G$ by $G_\varrho$, i.e., $G_\varrho = G\cap \Gamma_{n,\varrho}$ 
where $\Gamma_{n,\varrho}$ is 
$\ker \varphi_\varrho$ on $\GL(n,R)$.

There are two parts of our method, one for $\varphi_\varrho(G)$ and one 
for $G_\varrho\unlhd G$. Since $G_\varrho$ has finite index in a finitely 
generated group, it is finitely generated too. None of our 
algorithms call for a full generating set of $G_\varrho$
(which may be hard or even impossible to acquire), but rather a 
\emph{normal generating set}: a finite subset $N$ of $G_\varrho$ such 
that $G_\varrho$ is the normal closure $\langle N\rangle^G$ of 
$\langle N\rangle$ in $G$. There is a standard CGT procedure to obtain $N$; 
see \cite[pp.~299--300]{HoltBook}.
For this we need a presentation of $\varphi_\varrho(G)$ in the form
$\langle \varphi_\varrho(g_1), \ldots, \varphi_\varrho(g_r) \, |\,  
w_1, \ldots , w_l\rangle$ where $\{g_1, \ldots , g_r\}$ is the input 
generating set $S$ of $G$ and the $w_j$ are words in the 
$\varphi_\varrho(g_i)$. Efficient algorithms 
to compute such a presentation are available~\cite{Baaetal}. 
Replacing $\varphi_\varrho(g_i)$ by $g_i$ in each $w_j$, 
we get words $\tilde{w}_j$ over $S$.
Then $N = \{ \tilde{w}_1, \ldots , \tilde{w}_l\}$. 
Label this process  ${\tt NormalGenerators}(S,\varphi_\varrho)$.

We summarize our computational version of finite approximation.
\begin{enumerate}
\item Select a maximal ideal $\varrho$ of $R$.
\item Construct the congruence image $\varphi_\varrho(G)$
over the finite field $R/\varrho$. 
\item \label{StThree} Find a presentation of 
$\varphi_\varrho(G)$.
\item 
Compute a normal generating set of $G_\varrho$.
\end{enumerate}
The bulk of the computation is in step~\ref{StThree}; but 
it is done over a finite rather than infinite domain. 
This eases complexity issues such as matrix entry growth. 
We also gain access to the powerful algorithms for matrix 
groups over finite fields.

The clarity of our method flows from deep results 
such as Theorem~\ref{Malcev}, which ensure that the method
may be converted into an efficient algorithm.  
The solution of the orbit-stabilizer problem for nilpotent-by-finite groups
over $\Q$ presented in \cite{Dixon85} is a model for this kind of
computing with infinite linear groups. As noted previously, 
Dixon's paper was an impetus for subsequent work by Assmann, Eick,
and Ostheimer on computing with polycyclic groups over $\Q$. Other 
algorithms for infinite matrix 
groups that use congruence homomorphisms appear in \cite{Beals1,Beals2,RTB}. 
Those algorithms feature randomization, which we have so far eschewed.

\section{Deciding finiteness}\label{FinitenessSection}

An obvious launching point for 
the investigation of a potentially infinite group is to decide
whether or not it is finite.
Algorithms to test finiteness of finitely generated matrix
groups over $\Q$ were developed  in \cite{MR1173878,BBR}. These mix
deterministic and randomized techniques, and have 
integrality testing as a subprocedure. Finiteness testing of groups 
over fields other than $\Q$ is considered in 
\cite{Detinko01,Ivanyos01,RTB}; none of the algorithms
from those papers were implemented.
Ideas from \cite{BBR} are utilized in the {\sf GAP} package 
{\sf GRIM}~\cite{GRIM}. Both {\sf GAP} and {\sc Magma} use 
\cite{BBR} for their  default procedures to test finiteness over $\Q$.
In this section we present an algorithm to test finiteness over any field.

\subsection{Selberg--Wehrfritz theorems}\label{SWSubsection}

By a result of Schur's~\cite[p.~181]{SuprunenkoI}, if  the 
finitely generated subgroup $G$ of $\glnf$ is periodic (all elements
are torsion, i.e., have finite order) then $G$ is finite.
Usually we expect the torsion part of $G$ to be `small', which in 
characteristic zero means that $G$ is virtually torsion-free. 

We define more terminology. An element $g$ of $\glnf$ is 
\emph{unipotent} if it is conjugate to a unitriangular 
matrix in $\glnf$; equivalently, $g$ has characteristic polynomial 
$({\rm x}-1_n)^n$. If $\mathrm{char} \, \F = p>0$ then the 
unipotent elements of $\glnf$ are 
precisely its $p$-elements.
A subgroup $H$ of $\glnf$ is called unipotent if every element of $H$ 
is unipotent;  equivalently, $H$ may be conjugated into the group 
$\UT(n,\F)$ of all $n\times n$ upper unitriangular matrices over $\F$.
In positive characteristic $p$, the unipotent subgroups of $\glnf$ are 
the $p$-subgroups. Unipotent groups are nilpotent.

We now state a well-known key result for finitely generated 
linear groups (see \cite[Corollary~4.8, p.~56]{WBook}).
\begin{theorem}
\label{SWTheorem}
$G$ has a normal subgroup $H$ of finite index whose torsion elements are 
unipotent. In particular, if $\mathrm{char} \, \F = 0$ then $H$ is 
torsion-free.
\end{theorem}

Selberg proved Theorem~\ref{SWTheorem} for zero characteristic;
Wehrfritz extended it to all characteristics. The proof 
in \cite[p.~56]{WBook} 
does not give $H$ as a congruence subgroup; unlike the following 
(which implies Theorem~\ref{SWTheorem}).
\begin{proposition}[{\cite[Proposition~2.1]{Recognition}}] 
\label{Noetherian}
Let $\Delta$ be a Noetherian integral domain, and $\rho$ be a
maximal ideal of $\Delta$. If $g\in \Gamma_{n,\rho}$ has finite order
then $|g|$ is a power of $\mathrm{char} (\Delta/\rho)$.
\end{proposition}

Finitely generated integral domains are Noetherian. 
Proposition~\ref{Noetherian} specifies ideals 
$\varrho$ such that $G_\varrho=H$ 
as per Theorem~\ref{SWTheorem}. For such $\varrho$ we call $\varphi_\varrho$ 
an \emph{SW-homomorphism}. If $\mathrm{char} \, \F>0$ and 
$\varrho$ is any maximal ideal of $R$ 
then $\varphi_\varrho$ is an SW-homomorphism.

Proposition~\ref{Noetherian} is still not enough for our purposes. 
We also need
\begin{proposition}[{\cite[Theorem 4, p.~70]{SuprunenkoI}}] 
\label{RhoSquared}
Suppose that $\Delta$ is a Dedekind domain of characteristic zero, and 
$\rho$ is a maximal ideal of $\Delta$ such that 
$\mathrm{char}(\Delta/\varrho)$ is an odd prime $p$. If 
$p\not \in \varrho^2$ then $\Gamma_{n,\rho}$ is torsion-free.
\end{proposition}

If $\Delta = \Z$ then Proposition~\ref{RhoSquared} is a result of
 Minkowski: $\Gamma_{n,p}$ is torsion-free for odd primes $p$ (indeed 
$\Gamma_{n,m}$ is torsion-free for any odd integer 
$m$~\cite[Theorem~IX.8]{NewmanIntegral}).
\subsection{Constructing SW-homomorphisms}\label{ConstructSW}
We adhere to the notation and conventions of 
Subsection~\ref{ConstrCH}. 
\subsubsection{}
\label{RationalsSubsection}
If $R=\frac{1}{\mu}\mathbb Z$
then $\varphi_{1,p}$
is an SW-homomorphism for any odd prime $p$ not dividing $\mu$.

\subsubsection{}
\label{NumberFieldsSubsection}
Let $\mathbb{F} = \mathbb{P} = \mathbb{Q}(\alpha)$ for an algebraic 
number $\alpha$ with minimal polynomial $f(t)$ of degree $k$; then
$R \subseteq \frac{1}{\mu}\Z[\alpha ]$.
Let $p$ be a prime not dividing $\mu$, such that
either $p>nk+1$ (so that $\GL(n,\F)$ does not contain non-trivial 
$p$-subgroups) or $p$ does not divide the discriminant of $f(t)$.
It follows from Propositions~\ref{Noetherian} and \ref{RhoSquared}
that $\varphi_{2,p}$ is an SW-homomorphism
(see  \cite[p.~103]{Recognition}).

\subsubsection{}
\label{FunctionFieldsSubsection}

Let $\F = \E({\rm x}_1, \ldots, {\rm x}_m)$ where $\E$ is
$\Q$, $\P$, or $\F_q$. By Proposition~\ref{Noetherian},
$\varphi_{3,\alpha}$ and $\varphi_{3,\alpha,p}$
are SW-homomorphisms.

\subsubsection{}
\label{AlgebraicFunctionFieldsSubsection}

Let $\F = \eL(\beta)$ where $\eL = \E({\rm x}_1, \ldots, {\rm x}_m)$ 
and $|\F:\eL|=e$, so $R\subseteq \abk \frac{1}{\mu}\eL_0[\beta]$.
If $\mathrm{char}\, \F = 0$ then $\varphi_{4,\alpha}$ has 
torsion-free kernel by Proposition~\ref{Noetherian};
hence the composition $\varphi_{4,\alpha,p}$ of $\varphi_{4,\alpha}$ 
with $\varphi_p$ as in \ref{RationalsSubsection} or 
\ref{NumberFieldsSubsection} 
is an SW-homomorphism. 
If $\E=\F_q$
then $\varphi_{4,\alpha}$ is already an SW-homomorphism.

\bigskip

Thus, SW-homomorphisms always exist, and there are infinitely 
many $\varrho$ such that $\varphi_\varrho$ is an SW-homomorphism.
Moreover, if $\F$ is $\Q$ or $\P$ then $\varphi_\varrho$ 
is an SW-homomorphism for all but a finite number of $\varrho$ 
(cf.~\cite[Section~3.5]{Recognition}).

\subsection{The algorithm}\label{AFinTest}
The next lemma recaps the definition of 
SW-homomorphism.

\begin{lemma}\label{JustifyIsFinite}
Let $\varphi_\varrho$ be an SW-homomorphism on $G\leq \GL(n,R)$.
\begin{itemize}
\item[{\rm (i)}] Suppose that $\mathrm{char} \, R = 0$. Then $G$ is finite 
if and only if $G_\varrho =1$.
\item[{\rm (ii)}] Suppose that $\mathrm{char} \, R = p>0$. Then $G$ is finite 
if and only if $G_\varrho$ is a finite $p$-group
(i.e., is unipotent).
\end{itemize}
\end{lemma}

Lemma~\ref{JustifyIsFinite} guarantees correctness of the following.

\vspace*{14pt}

${\tt IsFinite}\hspace{.25pt} (S)$

\medskip

Input: a finite subset $S$ of $\GL(n,R)$, $\mathrm{char} \, R = p\geq 0$.

Output: ${\tt true}$ if $G=\gpess$ is finite; ${\tt false}$
otherwise.

\vspace*{1mm}

\begin{enumerate}

\item Select an ideal $\varrho$ of $R$ such that 
 $\varphi_\varrho$ is an SW-homomorphism, and construct $\varphi_\varrho(G)\leq 
\GL(n,q)$ where $|R/\varrho|=q$.

\item $N := {\tt NormalGenerators}\hspace{.25pt} (S, \varphi_\varrho)$.

\item \label{IsUSt3}
If $p=0$ and $N= 1$, or $p>0$ and $\langle N\rangle^G$ is unipotent,
then return ${\tt true}$;

\noindent else return ${\tt false}$.
\end{enumerate}

\vspace*{12pt}

If $\mathrm{char} \, R = 0$ then to confirm finiteness of 
$G$ it is enough to verify that $N$ (i.e., $G_{\varrho}$) is trivial.
If $p > 0$ then
we use the procedure ${\tt IsUnipotentClosure}(T,G)$ from 
\cite[Section~5.2]{Tits} in Step~\ref{IsUSt3}.
For a finite set $T\subseteq \glnf$, this tests 
whether $\langle T\rangle^G$ is unipotent. The computation is done in 
the enveloping algebra of $\langle T\rangle^G$ (smallest subalgebra
of $\mathrm{Mat}(n,\F)$ containing $\langle T\rangle^G$).
Testing unipotency of finitely 
generated subgroups of $\GL(n,\F)$ is easier; see, e.g., 
\cite[p.~108]{Large}.

As an auxiliary step, we might check whether randomly chosen elements 
(words over $S$) have finite order: by Schur's result,  
infinite $G$ contains elements of infinite order.
This has turned out to be a reliable way of certifying infiniteness
quickly.

For certain input $\tt IsFinite$ may be further modified.
If $\F$ is a function field then we only need the 
substitution homomorphism $\varphi_{3,\alpha}$ and computation with 
enveloping algebras, in place of ${\tt NormalGenerators}$. 
This may be helpful insofar as the words that arise are shorter than 
the words over $S$ that can arise in a run of ${\tt NormalGenerators}$.
See \cite{FinitenessZero,FinitenessPositive}.

Apart from being a practical algorithm valid over an arbitrary 
infinite field, ${\tt IsFinite}$ justifies decidability of the 
finiteness problem in the class of finitely generated linear groups.

\subsection{Recognition of finite matrix groups}

Let $G\leq \GL(n,\F)$ be finite.
If $\mathrm{char} \, \F = 0$ then an
SW-homomorphism $\varphi_\varrho$ maps $G$ isomorphically onto 
$\varphi_\varrho(G)$, a matrix group over some finite field $\F_q$.
If $\mathrm{char}\, \F =p>0$ then $G_\varrho$ could be a non-trivial
finite $p$-group. However, an
SW-isomorphism in this case may be obtained from an ad hoc recursion, 
as in \cite[Section~4.3]{Recognition}.

Once we have an isomorphic copy $\varphi_\varrho(G)\leq \GL(n,q)$ 
of $G$ in some $\GL(n,q)$, we `recognize' $G$ by subjecting 
$\varphi_\varrho(G)$ to the gamut of algorithms for matrix groups
over finite fields~\cite{Baaetal,BBS,gaprecog}. Amongst other things,
we can: compute $|G|$; test whether $G$ is solvable or nilpotent; 
compute a composition series of $G$; find a presentation of $G$; test 
membership of $g\in \glnf$ in $G$.

\section{Computing with virtually solvable linear groups}
\label{ComputingWithSF}

We move on to the next phase of 
investigating a finitely generated linear group.

\subsection{The Tits alternative}

Tits proved the following milestone result~\cite{Tits72}.
\begin{theorem}[The Tits alternative]\label{TAlt}
A finitely generated subgroup of $\GL(n,\F)$
either is solvable-by-finite, or contains a non-abelian free subgroup.
\end{theorem}
If $\mathrm{char}\, \F = 0$ then the conclusion of Theorem~\ref{TAlt} 
holds for all $G\leq \GL(n,\F)$.
Results analogous to the Tits alternative for 
other kinds of groups are given in \cite{BreuillardTits}, 
\cite[Section~2]{DixonTits}, \cite{LarsenShalev}, and
\cite[Section 4.5, pp.~154--162]{ShirvaniWehrfritz}.

The Tits alternative divides all finitely generated linear groups 
into two disparate classes, and it is vital that we are able to 
determine the class to which a given group belongs. An algorithm 
for doing this is given later in the section.

\begin{proposition}[{\cite[Corollary~10.18, p.~146]{WBook}}]
\label{WProp9} 
Linear groups satisfying the maximal condition on subgroups 
are polycyclic-by-finite, and vice versa.
\end{proposition}

Proposition~\ref{WProp9} points to the computational 
tractability of polycyclic-by-finite groups. 
It implies termination of a 
procedure to compute the normal closure $\langle N\rangle^G$ of a 
finite subset $N$ in a polycyclic-by-finite (linear) group $G$.

We note another condition for virtual solvability.
\begin{theorem}[{\cite[Theorem~10.9, p.~141]{WBook}}]
\label{WThm10} 
Suppose that each finitely generated subgroup of the linear 
group $G$ can be generated by $d$ elements, for some fixed
positive integer $d$. Then $G$ is solvable-by-finite.
\end{theorem}

Existing proofs of Theorem~\ref{TAlt} (as in, e.g., 
 \cite{DixonTits,Tits72})
do not translate into a practical algorithm to test virtual 
solvability (over any field $\F$). We proceed instead by way of
computational finite approximation. 

\subsection{Solvable-by-finite linear groups and congruence homomorphisms}
\label{StarredSubsection}
 
A \emph{block (upper) triangular group} in $\GL(n,\F)$ is a subgroup 
of the form
\[
H= \left(\renewcommand{\arraycolsep}{.1cm}
\begin{array}{cccc} H_1 & * & \cdots & *
\\
0 & H_2 & \cdots & *
\\
\vdots & \vdots & \ddots & \vdots\\
0 & 0 & \cdots & H_k \end{array} \right) 
\]
\noindent with zeros beneath the `diagonal part' 
$\mathrm{diag}(H_1, \ldots , H_k)$,
where $H_i\leq \abk \GL(n_i,\F)$.
If all $H_i=1$ then $H\leq \UT(n,\F)$; and if all
$n_i=1$ then $H\leq\mathrm{T}(n,\F)$, the
group of upper triangular matrices.
Any subgroup of $\GL(n,\F)$ is conjugate to one of the
above form 
where each $H_i$ is irreducible as a subgroup of $\GL(n_i,\F)$. 
If $H$ is unipotent-by-abelian then it is conjugate 
in $\mathrm{GL}(n,\overline{\F})$ to a subgroup 
of $\mathrm{T}(n,\overline{\F})$.
The set of all upper unitriangular matrices in $H$ is
a unipotent normal subgroup $U(H)$; 
this is the kernel of the projection of $H$ onto its 
`completely reducible part' $\mathrm{diag}(H_1, \ldots , H_k)$.
\begin{theorem}[Lie--Kolchin--Mal'cev]
\label{LKMTheorem}
Each solvable-by-finite linear group contains a 
unipotent-by-abelian subgroup 
of finite index.
\end{theorem}

There are algebraic and topological proofs of Theorem~\ref{LKMTheorem}
(see \cite[Theorem~7, p.~135]{SuprunenkoI} and 
\cite[Theorem~5.8, p.~77]{WBook}). 
The theorem implies that a solvable-by-finite
linear group can be conjugated
to block triangular form, 
with (irreducible) abelian-by-finite diagonal blocks.

The next result enables us to compute with a finite index 
unipotent-by-abelian subgroup 
of a given solvable-by-finite linear group.
\begin{theorem}[Wehrfritz~\cite{WJAlgebra}]
\label{WJAThm}
Let $G \leq \GL(n,R)$ be solvable-by-finite,
 and let $\varrho$ be an ideal of $R$. 
Then $G_\varrho$ is unipotent-by-abelian if 
\begin{itemize}
\item[{\rm (i)}] 
$R/\varrho$ has prime characteristic greater than $n$; or
\item[{\rm (ii)}] $R$ is a Dedekind domain of characteristic zero, 
$\varrho$ is a maximal ideal of $R$, $R/\varrho$ has odd 
characteristic $p$, and $p\not \in \varrho^{p-1}$. 
\end{itemize}
\end{theorem}

\begin{remark}
 $G_\varrho$ in Theorem~\ref{WJAThm}~(ii) is Zariski-connected 
(see Subsection~\ref{ZDAndAG}).
\end{remark}
 
Theorem~\ref{WJAThm} is proved in \cite[Lemma 9]{Dixon85} for 
$\F = \Q$.
This was background for 
the algorithm in \cite{BettinaBjorn} to test virtual solvability over $\Q$.
The Monte Carlo algorithm of \cite{Beals1} to decide the Tits alternative 
over $\Q$ relies on solvability testing of matrix groups over 
finite fields~\cite{Luks}. 

\subsection{A computational version of the Tits alternative}
\label{TAltRealization}

Call $\varphi_\varrho$ for an ideal $\varrho\subseteq R$ as
in Theorem~\ref{WJAThm} a \emph{W-homomorphism}.

\subsubsection{}
\label{WOne}
$\varphi_{1,p}$ for an odd prime $p \in \mathbb Z$ 
as in \ref{RationalsSubsection} is a W-homomorphism.
\subsubsection{}
\label{WTwo}
$\varphi_{2,p}$ 
as in \ref{NumberFieldsSubsection} is a W-homomorphism
if either $p>n$, or $p$ is coprime to the discriminant of 
the minimal polynomial $f(t)$ of $\alpha$.
\subsubsection{}
If $\mathrm{char}\, \F=0$ then  
$\varphi_{3,\alpha,p}$ is a W-homomorphism;
if $\mathrm{char}\, \F=p>n$ then the substitution map 
$\varphi_{3,\alpha}$ on its own is
a W-homomorphism (see \ref{FunctionFieldsSubsection}).
\subsubsection{}
\label{WFour}
If $\mathrm{char}\, \F = 0$ then $\varphi_{4,\alpha,p}$
is a W-homomorphism; if $\E=\F_q$ then $\varphi_{4,\alpha}$ is a 
W-homomorphism
(see \ref{AlgebraicFunctionFieldsSubsection}).

\bigskip

Our computational realization of the Tits alternative follows.

\vspace{10pt}

${\tt IsSolvableByFinite}(S)$

\medskip

Input: $S = \{g_1, \ldots, g_r\}\subseteq \mathrm{GL}(n,R)$.

Output: ${\tt true}$ if $G = \gp{S}$ is solvable-by-finite;
${\tt false}$ otherwise.

\vspace{1mm}

\begin{enumerate}

\item Select $\varrho\subseteq R$ such that $\varphi_\varrho$ is a
W-homomorphism, and construct $\varphi_\varrho(G)$.

\item \label{SecondStep} $N := {\tt NormalGenerators}(S,\varphi_\varrho)$.

\item \label{ThirdStep} Return $\tt true$ if $\langle N\rangle^G$ is 
unipotent-by-abelian; else return $\tt false$.

\end{enumerate}

\vspace{12pt} 

Step~\ref{ThirdStep} is a matrix algebra computation.
(The normal closure cannot be computed directly
by a standard recursion, as this may not terminate 
if $G$ is not polycyclic-by-finite; cf.~Proposition~\ref{WProp9}.)
To find a basis of the enveloping algebra $\langle G_\varrho\rangle_\F$, 
or to test whether $\langle N\rangle^G$ is unipotent-by-abelian,
only a normal generating set for $G_\varrho$ is required. We
already have this from Step~\ref{SecondStep}. 
Even if we can find a full generating set of $G_\varrho$, 
the enveloping algebra method may still be preferable to a
direct normal closure computation (see the
penultimate paragraph of  Subsection~\ref{AFinTest}).    

\subsection{Other group-theoretic properties}
\label{FourFour}

We now test narrower attributes of the 
input solvable-by-finite linear group $G$: whether it is  
nilpotent-by-finite, abelian-by-finite, central-by-finite,
solvable, nilpotent. 
Maintaining a common theme, we give practical 
algorithms that justify decidability.

A class ostensibly not too far removed from finitely generated 
solvable-by-finite groups is polycyclic-by-finite groups. 
However, there is in fact a large distance between the classes: 
a polycyclic-by-finite group is finitely presentable, has every 
subgroup finitely generated, and satisfies the maximal condition 
on subgroups; whereas none of this is true in general for 
solvable-by-finite groups. So algorithmic methods for 
polycyclic-by-finite groups may not work at all for 
solvable-by-finite groups (cf.~the comments after 
Proposition~\ref{WProp9}).

To test solvability of $G$ we add checking solvability of
$\varphi_\varrho(G)$ to $\tt IsSolvable$-$\tt ByFinite$. This 
single extra step
is readily accomplished; see~\cite{BettinaBjorn,Baaetal,Luks}.

For nilpotent-by-finite groups, we combine some theory 
of nilpotent linear groups with the next result.
\begin{proposition}[{\cite[Corollary~5.2]{Tits}}]
\label{LemCor52Tits}
Suppose that $R$ is a Dedekind domain of characteristic zero, and
$\varrho$ is a maximal ideal of $R$ such that
$\mathrm{char}(R/\varrho)=p > 2$, where $p \notin \varrho^{p-1}$. Then
$G \leq \mathrm{GL}(n, R)$ is nilpotent-by-finite (respectively,
abelian-by-finite) if and only if $G_{\varrho}$ is nilpotent (respectively,
abelian).
\end{proposition}

Proposition~\ref{LemCor52Tits} 
follows from Theorem~\ref{WJAThm}~(ii). Since it is connected,
$G_\varrho$ is nilpotent (respectively, abelian)
if it is nilpotent-by-finite (respectively, abelian-by-finite).

Denote by $g_d$, $g_u\in \mathrm{GL}(n,\F)$ the diagonalizable and
unipotent parts of $g\in \glnf$. So $g_d$ is conjugate to a diagonal 
matrix over $\overline{\F}$, $g_u$ is unipotent,
and $g =g_dg_u=g_ug_d$.
This is the Jordan decomposition of $g$;
see \cite[Theorem~7.2, p.~91]{WBook}.
For $X\subseteq \glnf$ 
set $X_d= \{ h_d\mid h\in G\}$,
$X_u=\abk \{ h_u \mid h\in G\}$.
\begin{lemma}\label{PropertiesOfNilpotent} 
Let $G=\langle S\rangle\leq \glnf$.
\begin{itemize}
\item[{\rm (i)}] $G$ is nilpotent if and only if 
$\langle S_d\rangle$, 
$\langle S_u\rangle$ are nilpotent 
and centralize each other.
\item[{\rm (ii)}] If $G$ is nilpotent then $G\leq \langle S_d\rangle
\times \langle S_u\rangle$.
\end{itemize} 
\end{lemma}
Proofs of Lemma~\ref{PropertiesOfNilpotent} are
in \cite{Nilpotent} and \cite{SuprunenkoI}.

Now we can state an algorithm to test virtual nilpotency, based
on the above.

\vspace{12pt}

${\tt IsNilpotentByFinite}(S)$

\medskip

Input: a finite subset $S$ of $\mathrm{GL}(n,R)$, $R$ a Dedekind
domain of characteristic $0$.

Output: ${\tt true}$ if $G = \gp{S}$ is nilpotent-by-finite; otherwise
${\tt false}$.

\medskip

\begin{enumerate}

\item \label{PreCrunch}
Select $\varrho$ such that $\varphi_\varrho$ is a W-homomorphism.

\item $N := {\tt NormalGenerators}(S,\varphi_\varrho)$.

\item \label{Crunch} If 
$\langle x^g : x\in N_d, \, g\in G \rangle$ is abelian, 
$\langle y^g : y\in N_u, \, g\in G \rangle$is unipotent,  
and these two groups commute elementwise,
then return $\tt true$; else return
$\tt false$. 

\end{enumerate}

\vspace{12pt}

Step~\ref{Crunch}  uses  
$\tt IsAbelianClosure$ and
${\tt IsUnipotentClosure}$ from \cite[p.~404]{Tits}.
Once more these involve computation in related 
enveloping algebras.

The procedure ${\tt IsAbelianByFinite}(S)$ tests whether 
$G=\langle S\rangle\leq \GL(n,R)$ is abelian-by-finite, where 
$R$ is a Dedekind domain of characteristic zero. 
All steps are the same as in 
${\tt IsNilpotentByFinite}$ except for Step~\ref{Crunch}, which now 
simply returns ${\tt IsAbelianClosure}(N, S)$.

Next we show how to decide whether $G$ is central-by-finite.
\begin{lemma}[{\cite[Corollary~5.8]{Tits}}]
\label{Transfer}
Let $G\leq \GL(n,R)$ where $\mathrm{char} \, R=0$, and let 
$\varphi_\varrho$ be an SW-homomorphism on $\GL(n,R)$. 
Then $G$ is central-by-finite if and only if $G_\varrho$ is 
central in $G$.  
\end{lemma}
If $G$ is central-by-finite then the commutator subgroup $[G,G]$
generated by all $[x,y]= x^{-1}y^{-1}xy$ is finite.  
The non-trivial direction of Lemma~\ref{Transfer} 
follows from this and $G_\varrho$ being a torsion-free 
normal subgroup of $G$. Thus ${\tt IsCentralByFinite}$ returns 
${\tt true}$ if the input $S$ centralizes 
${\tt NormalGenerators}(S,\varphi_\varrho)$ 
and $\tt false$ otherwise.

We round out the section with nilpotency testing
(in characteristic zero).
This will not be a simple modification of
$\tt IsNilpotentByFinite$, as an extension of one nilpotent 
group by another need not be nilpotent (cf.~testing solvability via
$\tt IsSolvableByFinite$).
\begin{lemma}[{\cite[Lemma~4.9]{Nilpotent}}]
\label{CRNilpotentIsCF}
Let $G\leq \GL(n,R)$ be nilpotent, $\mathrm{char} \, R=0$, and 
suppose that $G=G_d$. If $\varphi_\varrho$ is an SW-homomorphism 
then $G_\varrho\leq Z(G)$.
\end{lemma}

Lemmas~\ref{PropertiesOfNilpotent} and \ref{CRNilpotentIsCF}
lead to the following.
 
\vspace{12pt}

${\tt IsNilpotent}(S)$

\medskip

Input: a finite subset $S$ of $\mathrm{GL}(n,\F)$, $\mathrm{char}\, \F = 0$.

Output: ${\tt true}$ if $G = \gp{S}$ is nilpotent; otherwise
${\tt false}$.

\medskip

\begin{enumerate}

\item $H:=\langle S_d\rangle$, $K:=\langle S_u\rangle$.

\item If $K$ is not unipotent, or $[H,K]\neq 1$, then 
return $\tt false$.

\item \label{TestOverFiniteField}
Select $\varrho$ such that
 $\varphi_\varrho$ is an SW-homomorphism on $\GL(n,R)$.

\noindent If $\varphi_\varrho(G)$ is not nilpotent then return $\tt false$.

\item If $H_{\varrho}\not \leq Z(H)$ then return $\tt false$;
else return $\tt true$.

\end{enumerate}

\vspace{12pt}

For nilpotency testing over finite fields 
(Step~\ref{TestOverFiniteField}), see \cite{Large}. 
The papers \cite{Large,Nilpotent}
contain many more algorithms for nilpotent linear groups. 

\subsection{Structure of solvable-by-finite linear groups}
\label{InStructureSF}
 
Further study of a solvable-by-finite subgroup 
$G=\langle S\rangle$ of $\glnf$ 
begins by computing its main structural components.  
This leads inevitably to a consideration of ranks.

We may assume that $G$ is block triangular 
with completely reducible abelian-by-finite diagonal part 
(cf.~Subsection~\ref{StarredSubsection}). 
Let $\pi$ be the projection of $G$ onto its diagonal part.
Then $\ker \, \pi=U(G)$, the unipotent radical of $G$.
Certainly $\pi(G)$ is finitely generated, whereas $U(G)$ is
finitely generated if and only if $G$ is polycyclic-by-finite.

We note two related procedures: $\tt IsCR$ 
and $\tt CRPart$. $\tt IsCR$ tests whether $G$ is completely 
reducible; equivalently, whether $U(G)=\abk 1$. 
$\tt CRPart$ returns a generating set of $\pi(G)$. Here only 
the least complicated case $\mathrm{char}\, R = 0$ is
reviewed. 

Let $\varphi_\varrho$ be a W-homomorphism for $G$. Then $G$ is 
completely reducible if and only if $G_\varrho$ is completely 
reducible abelian. The latter can be 
tested by customary manipulations in an enveloping algebra  
of $N= {\tt NormalGenerators}(S,\varphi_\varrho)$; 
see \cite[Section~4]{Tits}.
When $G$ is nilpotent-by-finite, $G$ is completely reducible if and 
only if $N_u =1$;  when $G$ is nilpotent, it suffices to 
check whether $S_u=1$.
$\tt CRPart$ is described in \cite[Section~4.2]{FiniteRank}.
The  completely reducible 
part of nilpotent $G$ is $\langle S_d\rangle$.

So we can decide whether $G$ is completely reducible. 
We reiterate that although $U(G)$ may 
not be finitely generated, it is 
nilpotent, and torsion-free in characteristic zero. 
Consequently $U(G)$ has finite rank (see below), and
some headway could be made computationally 
using P.~Hall's methods for infinite nilpotent 
groups (cf.~\cite[pp.~30--33]{LennoxRobinson}).

\subsection{Finite rank linear groups}
\label{NewFRSection}

This subsection illustrates how rank restrictions facilitate 
computing with finitely generated solvable-by-finite 
linear groups. 

Recall that a group $H$ has \emph{finite Pr\"{u}fer rank} $\rk(H)$ if
each finitely generated subgroup  can be generated by
$\rk(H)$ elements, and $\rk(H)$ is the least such integer.
Finite rank linear groups are solvable-by-finite (Theorem~\ref{WThm10}). 
The converse is false.
\begin{example}\label{Example413}
\emph{The subgroup 
$\big \langle \hspace{.1pt}
\mathrm{diag}(1,{\rm x}),
{\tiny 
\renewcommand{\arraycolsep}{.1cm}
\left(\begin{array}{cc}
1 & 1  \\
0 & 1 
\end{array} \right)} 
\big \rangle$ of $\GL(2,\F_2({\rm x}))$ is solvable,
but does not have finite Pr\"{u}fer rank, as it contains
${\tiny \renewcommand{\arraycolsep}{.1cm}
\left(\begin{array}{cc}
1 & {\rm x}^k  \\
0 & 1 
\end{array} \right)}$ for all $k\geq 1$.}
\end{example}
\begin{proposition}[{\cite[Corollary~2.5]{FiniteRank}}]
\label{SFSufficesForQLinear}
A finitely generated subgroup of $\glnf$ 
has finite Pr\"{u}fer 
rank if and only if it is solvable-by-finite 
and $\Q$-linear, i.e., isomorphic to 
a subgroup of $\GL(m,\Q)$ for some $m$.
\end{proposition}
\begin{example}
\emph{Polycyclic-by-finite groups are $\Z$-linear,
 so have finite Pr\"{u}fer rank.}
\end{example}

A group $H$ has \emph{finite torsion-free rank} if there is 
a subnormal series of finite length in $H$, with
each factor either periodic or infinite cyclic. 
The number $\hi(H)$ of infinite cyclic factors 
is the {\em torsion-free rank} or {\em Hirsch number} 
of $H$.
\begin{proposition}[{\cite[Proposition~2.6]{FiniteRank}}]
For a finitely generated subgroup $G$ of $\GL(n,\Q)$ the
following are equivalent.
\begin{itemize}
\item[{\rm (i)}] $G$ is solvable-by-finite.
\item[{\rm (ii)}] $G$ has finite Pr\"{u}fer rank.
\item[{\rm (iii)}] $G$ has finite torsion-free rank.
\end{itemize}
\end{proposition}
\begin{proposition}[Cf.~{\cite[Sections~3.2 and 4.3]{FiniteRank}}]
\label{Isolator}
Each finitely generated solvable-by-finite 
subgroup $G$ of $\GL(n,\Q)$ has a
finitely generated subgroup $H\leq U(G)$ such that 
$U(G)=H^G$ and $\hi(H)=\hi(U(G))=\rk(U(G))$.
\end{proposition}
In Proposition~\ref{Isolator}, $U(G)$ is the \emph{isolator} of $H$
in $U(G)$: for each $g\in U(G)$ there is
a positive integer $m$ such that $g^m\in H$
(see \cite[Section~2.1]{LennoxRobinson}).

Let $\P$ be a number field. 
To test whether $G=\langle S\rangle\leq \GL(n,\P)$ has 
finite Pr\"{u}fer rank, we just run $\tt IsSolvableByFinite$.
Suppose that $G$ is  solvable-by-finite. 
If $G$ is completely reducible then it is  abelian-by-finite, 
$G_\varrho$ is (torsion-free) abelian, and 
$\hi(G) = \hi(G_\varrho)$. The rank of 
$G_\varrho$ may be computed using algorithms 
to construct presentations of irreducible abelian subgroups of 
$\GL(n,\P)$~\cite[Subsections~4.1.1 and 4.5]{FiniteRank}. 
This gives a procedure ${\tt RankCR}$  that returns $\hi(G)$
for input completely reducible $G$.
We can find $\hi(G)$ with similar ease when $G$ is unipotent
(see \cite[Subsection~4.1.2]{FiniteRank});
call the associated procedure ${\tt RankU}$. 
Then ${\tt RankRadical}$ computes $\hi(U(G))=\rk(U(G))$ 
as follows. 
First, a presentation of ${\tt CRPart}(S)$ is employed
to produce a set $Y$ of normal generators for $U(G)$. Secondly, 
a generating set of $H$
as in Proposition~\ref{Isolator} is found
by the method in \cite[Section~4.3]{FiniteRank};
${\tt RankU}$ is required and the computation initializes at $Y$. 
Then $\hi(U(G))=\hi(H)$.
Furthermore, since $\hi(G) = \hi(G/U(G)) + \hi(U(G))$, and 
$G/U(G)$ is isomorphic to ${\tt CRPart}(S)$,
we obtain a procedure ${\tt HirschNumber}$ that accepts 
$S$ and returns $\hi(G)$.

The group $H\leq U(G)$ constructed in ${\tt RankRadical}$
deserves more scrutiny. 
It has a series $1= H_k \lhd H_{k-1} \lhd \cdots \abk \lhd H_0=H$ 
with infinite cyclic factors $H_{i-1}/H_i$. Let 
$H_{i-1} = \langle u_i,H_i\rangle$. 
Then since $U(G)$ 
is the isolator of $H$, we may represent each  
$g\in U(G)$ uniquely as a $k$-tuple of 
canonical rational parameters $(\alpha_1, \ldots , \alpha_k)$ 
where $g=u_1^{\alpha_1} \cdots u_k^{\alpha_k}$;
 see \cite[pp.~29--34]{LennoxRobinson}.
Thus $\{ u_1, \ldots , u_k\}$ serves as a `basis' of $U(G)$.

The next theorem expands on a result by D.~J.~S.~Robinson 
about finitely generated
solvable groups of finite abelian ranks.
\begin{theorem}[{\cite[Theorem~3.1]{FiniteRank}}]
\label{Thm3Point1}
Let $H\leq G\leq \glnf$ where $G$ is finitely generated and of
finite Pr\"{u}fer rank. Then $|G:H|<\infty$ if and only if 
$\hi(H)=\hi(G)$.
\end{theorem}

Theorem~\ref{Thm3Point1} yields the procedure ${\tt IsOfFiniteIndex}$.
This accepts generating sets $S_1$, $S_2$ 
for solvable-by-finite subgroups $G$, $H$ of $\GL(n,\P)$, 
respectively, such that 
$H\leq G$; and returns $\tt true$ if and only if
${\tt HirschNumber}(S_1)={\tt HirschNumber}(S_2)$.

The final decision problem that we discuss here has ties to 
Section~\ref{FinitenessSection} and 
the next section.

A subgroup of $\GL(n,\Q)$ is \emph{integral} if it  
can be conjugated into $\GL(n,\Z)$.  Finite subgroups of 
$\GL(n,\Q)$ are integral \cite[p. 46]{SuprunenkoI}.
The following integrality criterion is used in 
one of the finiteness testing algorithms 
from \cite{BBR} 
(valid for input 
over a quotient field $\F$ of a principal ideal domain).
\begin{proposition}\label{AnotherIntegralityCriterion}
Suppose that each element of $G\leq \GL(n,\Q)$ has 
trace in $\Z$. Then $G$ is integral if and only if either 
$G$ is finitely generated or the enveloping algebra
$\langle G \rangle_\Q$ is semisimple.
\end{proposition}

Integrality testing may be 
easier when the input is solvable-by-finite.
\begin{lemma}[{\cite[Lemma~4.1]{DeGrDetF}}]
\label{22ArithmSolvable}
Let $G\leq \GL(n,\Q)$ be  finitely generated 
solvable-by-finite,
and let $p$ be an odd prime such that 
$\varphi_p$ is a W-homomorphism on $G$.
Then $G$, $\pi(G)$, and $\pi(G_p)$ are all integral if any
one of these groups is integral.
\end{lemma}

Lemma~\ref{22ArithmSolvable} gives us the following
procedure from \cite[Section~4]{DeGrDetF}, 
applied there in arithmeticity 
testing (see Subsection~\ref{522}). 

\vspace{12pt}

${\tt IsIntegralSF} (S)$

\medskip

Input: a finite subset $S$ of $\GL(n,\Q)$ such that $G=\gpess$ is
solvable-by-finite.

Output: ${\tt true}$ if $G$ is integral; ${\tt false}$ otherwise.

\vspace{1mm}

\begin{enumerate}

\item  $N := {\tt NormalGenerators}(S,\varphi_p)$,
 $\varphi_p$ a W-homomorphism.

\item Return ${\tt
true}$ if each $g\in N$ is integral (the 
characteristic polynomial of $g$ has integer coefficients
and $\mathrm{det} (g) = \pm 1$);
else return ${\tt false}$.
\end{enumerate}

\vspace{6pt}

\begin{remark}
For finite $G$ we have $N=1$, and ${\tt IsIntegralSF}$
returns $\tt true$ as expected.
\end{remark}

\subsection{Implementation}

All algorithms from Sections~\ref{FinitenessSection} 
and \ref{ComputingWithSF} have been implemented 
in {\sc Magma} by Eamonn O'Brien and the authors.
Sample experiments are given in the documentation 
for \cite{Infinite}. 

\section{Computing with Zariski dense and arithmetic groups}
\label{ComputeDenseArith}

\subsection{Linear groups with a non-abelian free subgroup}
\label{LGWithFreeNonAbelian}

Until now we have been concerned almost exclusively 
with virtually solvable linear groups. 
Since they are built up from infinite abelian and 
finite blocks of restricted structure, intuition might 
suggest that these groups are rare.
Indeed, for various $\F$, it is unlikely
that a randomly selected finite subset of $\glnf$ generates a
solvable-by-finite group~\cite{Aoun,LarsenShalev,Rivin}. 

Linear groups that are not virtually solvable pose computational 
challenges of a different nature to those posed by virtually solvable groups. 
For example, although membership testing in a finitely generated 
solvable-by-finite subgroup of $\GL(n,\Q)$ is decidable~\cite{Kopytov2}, there 
exist subgroups of $\GL(4,\Z)$ in which 
membership testing is undecidable~\cite{M58}.
Hence this and related algorithmic problems in $\GL(n,\Z)$
for $n\geq 4$ (e.g., subgroup conjugacy testing) 
are  undecidable; see \cite[Section~5]{MillerSurvey}
and \cite[Section~3]{Survey}.
In Subsection~\ref{NonSFAndFurtherOP}, we note some other 
problems  in the class of non-solvable-by-finite linear groups
where decidability is unknown.

\subsection{Zariski density and arithmetic groups}\label{ZDAndAG}

To make computation with non-solvable-by-finite linear groups
feasible, we impose some natural
conditions on the input (which 
need not generate a non-solvable-by-finite group). 

First we give some definitions.
A subset $S$ of an $m$-dimensional $\F$-space $V$ is
 \emph{algebraic} if there exists nonempty 
$F\subseteq \mathbb{F}[{\rm x}_1, \ldots ,{\rm x}_m]$ 
such that $S$ is the set of zeros of all polynomials 
in $F$.
The Zariski topology on $V$ is the topology whose closed sets
are the algebraic subsets.
An \emph{algebraic group}  is a subgroup of
$\GL(n,\F)$ that is closed in 
the $n^2$-dimensional space $V=\mathrm{Mat}(n,\F)$.
The essential case for us is algebraic $\Q$-groups $\mathcal{G}\leq \GL(n,\C)$, 
i.e., $\mathcal G$ is the set of mutual zeros of a collection
of polynomials over $\Q$.

When $n=2s$ is even, the symplectic group
$\mathrm{Sp}(n,R)$ is defined to be 
\[
\{ h\in \GL(n,R) \, |\, hJh^\top = J\} \quad \text{where} 
\quad
J = {\footnotesize \begin{pmatrix}
 \phantom{-} 0_s & 1_s \\
-1_s & 0_s 
\end{pmatrix}}.
\]
Each of $\GL(n,\C)$, $\SL(n,\C)$, and $\Sp(n,\C)$ is an algebraic 
$\Q$-group.

Any linear group $G\leq \glnf$ is a subgroup of some linear algebraic 
group; say $\GL(n,\F)$ itself. 
The `smallest' algebraic group in $\glnf$ containing
$G$ is its \emph{Zariski closure} $\overline{G}$.
An algorithm to compute $\overline{G}$ for finitely generated 
$G$ and infinite $\F$ is given in \cite{Derksen}.
We assume that $G$ is a Zariski dense subgroup of an
algebraic group $\mathcal G$, i.e., $\mathcal{G} =\overline{G}$.
If we wish to compute with non-solvable-by-finite $G$ then we
may also assume that $\mathcal G$ is non-solvable. 

Now we introduce an important class of dense subgroups.
If $G\leq \GL(n,\C)$ and $R$ is a subring of $\C$ then
 $G_R:=G\cap \GL(n,R)$. 
Let $\mathcal G$ be an algebraic $\Q$-group.
We say that $H\leq \mathcal{G}_\Q$ is an \emph{arithmetic subgroup
of $\mathcal G$} (or merely \emph{arithmetic} when $\mathcal G$ is 
understood) if $H$ is commensurable with $\mathcal{G}_\Z$, meaning 
that
$|H:H_\Z|$ and  $|\mathcal{G}_{\Z}:H_\Z|$ are finite.
In particular, finite index subgroups of $\mathcal{G}_{\Z}$ are 
arithmetic. If $\mathcal{G}$ is semisimple then an arithmetic subgroup 
$H\leq \mathcal{G}$ is not solvable-by-finite, and we call $H$ 
a \emph{semisimple arithmetic group} \cite[p.~91]{LubotzkySegal}.

By a famous result of Borel and 
Harish-Chandra~\cite[p.~134]{ZalesskiiII},
arithmetic subgroups of algebraic $\Q$-groups are finitely 
presentable.
We remark that the notion of arithmetic group may be framed in algebraic
groups $\mathcal G$ over fields $\F$ other than $\Q$. 
If $\mathrm{char}\, \F\neq 0$ then an arithmetic subgroup need not be 
finitely presentable, nor even finitely 
generated; examples are $\SL(n, \F_q[{\rm x}])$ for $n= 3$, $2$ respectively
(see \cite[p.~2981]{ZalesskiiIII} and \cite[p.~134]{ZalesskiiII}).

Arithmetic subgroups in $\mathcal{G}_\Z$ are Zariski 
dense. A subgroup of
$\mathcal{G}_\Z$ that is dense but not arithmetic is 
a \emph{thin matrix group}~\cite{Sarnak}.
These are ubiquitous in $\SL(n,\Z)$~\cite{Fuchs,Rivin}.
 
See \cite[Chapters~5 and 14]{WBook}, 
\cite[Chapter~1, Sections~5 and 6]{ZalesskiiII} for more on 
algebraic groups, and \cite{Humphreys} for more on arithmetic 
groups.

\subsection{Decidability for arithmetic groups}\label{522}

Let $\mathcal G$ be an algebraic $\Q$-group. 
An arithmetic group $H\leq \calg_\Z$ 
is \textit{explicitly given} if  membership in $H$ of 
each $g\in \calg_{\Z}$ can be tested, and an upper bound on 
$|\calg_{\Z} : H|$ is known.
Grunewald and Segal~\cite{GSI,GSII,GSIII} 
justified decidability of problems for explicitly given
arithmetic groups. 
One of these is constructing
a (finite) generating set of an arithmetic subgroup of $\calg_{\Z}$. 
Note that the algorithm in \cite{GSI} to construct a generating 
set from input polynomials is not always practical.
Also, sometimes we will have a generating set of $\calg_\Z$ 
\emph{a priori}; 
e.g., if $\calg_\Z = \SL(n,\Z)$ or $\Sp(n,\Z)$.

Algorithms to construct a generating set of $\calg_{\Z}$
for unipotent or abelian $\calg$ are given 
in \cite{deGraaf1,deGraaf2}.
These are ingredients in the procedure 
${\tt GeneratingArithmetic}$ from 
\cite{DeGrDetF}, which constructs a generating set of a finite index 
subgroup of $\calg_{\Z}$ for a solvable algebraic $\Q$-group
$\calg$.
Motivation for ${\tt GeneratingArithmetic}$ is supplied by

\vspace{4pt}

\begin{quote}
(AT) \ \emph{Arithmeticity testing}: if $H$ is a
finitely generated subgroup of $\calg_{\Z}$, 
determine whether
$|\calg_{\Z}:H|$ is finite.
\end{quote}

\vspace{4pt}

\noindent 
It is  unknown whether (AT) is decidable for all $\calg$. 
However, since  $\calg_{\Z}$ is finitely presentable, if $|\calg_{\Z}:H|$ is 
finite then this can be detected by Todd--Coxeter coset enumeration (not an 
advisable option in practice;
Todd--Coxeter may not terminate even for very small input).
Thus (AT) is semidecidable as per \cite[p.~149]{HoltBook}.

When $\calg$ is solvable, arithmetic groups are polycyclic, 
and it was proved in \cite{DeGrDetF} that (AT) is decidable.
\begin{proposition}\label{Recapitulate}
Let $\calg$ be a solvable algebraic $\Q$-group.
A finitely generated subgroup $H$ of $\calg$ is arithmetic 
if and only if $H$ is integral and $\hi(H) = \hi(\calg_\Z)$. 
\end{proposition}
Therefore $H\leq \calg_\Z$ is arithmetic if and only if
$H$ and $\calg_\Z$ have the same Hirsch number 
(cf.~Theorem~\ref{Thm3Point1}).
 Proposition~\ref{Recapitulate} gives 

\vspace{12pt}

${\tt IsArithmeticSolvable}\hspace{.25pt} (S, \calg)$

\medskip

Input: a finite subset $S$ of $\calg_\Q$,  $\calg$ a solvable algebraic
$\Q$-group.

Output: ${\tt true}$ if $H=\gpess$ is arithmetic; ${\tt false}$
otherwise.

\vspace*{3mm}

\begin{enumerate}

\item If ${\tt IsIntegralSF}\hspace{.25pt} (S) = {\tt false}$ then
return ${\tt false}$.

\item $T := {\tt GeneratingArithmetic}\hspace{.25pt} (\calg)$.

\item If ${\tt HirschNumber}\hspace{.25pt} (S)\neq {\tt
HirschNumber}\hspace{.25pt} (T)$ then return ${\tt false}$;

else return ${\tt true}$.
\end{enumerate}

\vspace*{11pt}

To test arithmeticity when $\calg$ is unipotent we only need to
compare ${\tt Hirsch}$-${\tt Number}\hspace{.25pt} (S)$ and 
${\tt HirschNumber}\hspace{.25pt} (T)$, i.e., the first step could
be omitted.

\subsection{The congruence subgroup property}\label{CSPSection}

The class of arithmetic subgroups of algebraic $\Q$-groups---even 
just the semisimple ones---is very wide. A comparison of
$\SL(n,\Q)$ and $\SL(n,\Z)$ reveals how we 
should limit our scope. In $\SL(n,\Q)$ a proper normal subgroup
is scalar and has order at most $2$, 
whereas $\SL(n,\Z)$ has a plurality of normal subgroups, e.g., the 
(principal) congruence subgroups 
$\Gamma_{n,m}:= \ker \, \varphi_m$ for all positive integers $m$ 
(terminology and notation as in Subsection~\ref{FASubsection}).
The complete inverse image of the centre of $\SL(n,\Z/m\Z)$ in 
$\SL(n,\Z)$ under $\varphi_m$ is an example of a normal subgroup 
that is not a congruence subgroup.
The question of whether $\SL(n,\Z)$ has normal subgroups of infinite index, 
or more generally normal subgroups not containing any PCS, was raised 
as long ago as the 19th century. Fricke and Klein constructed subgroups 
of finite index in $\SL(2,\Z)$ that do not contain any PCS (see 
\cite{Raghunathan}). Also $\SL(2,\Z)$ has normal subgroups of infinite 
index: e.g., the normal closure of 
%
$\big\langle
{\tiny \renewcommand{\arraycolsep}{.1cm}
\left(\begin{array}{cc}
1 & m  \\
0 & 1 
\end{array} \right)}
\big\rangle$
for $m>5$ (see \cite[p.~33]{Mennicke}).
The story is very different
for degrees greater than $2$.
\begin{theorem}[Cf.~\cite{Bass}]
For $n>2$, each normal non-central subgroup of $\, \SL(n,\Z)$ or $\Sp(n,\Z)$ 
contains a PCS. 
\end{theorem}
\begin{corollary}\label{CSPDefSLSp}
Each finite index subgroup of $\SL(n,\Z)$ or $\mathrm{Sp}(n,\Z)$ for $n>2$ 
contains a PCS.
\end{corollary} 

Corollary~\ref{CSPDefSLSp} says that $\SL(n,\Z)$ or $\mathrm{Sp}(n,\Z)$ 
both have the \emph{congruence subgroup property} (CSP).
How prevalent is the CSP? We have seen that it does not hold for
$\SL(2,\Z)$. 
But if $n>2$ and $\P$ is a number field that is not totally imaginary, 
then a finite index subgroup of $\SL(n,\mathcal{O}_\P)$ or 
$\Sp(n,\mathcal{O}_\P)$ contains $\Gamma_{n,\varrho}$ for some 
maximal ideal $\varrho$ of $\mathcal{O}_\P$~\cite{Bass}. 
Determining whether arithmetic subgroups of an
algebraic group have the CSP is
`the congruence subgroup problem' (see \cite{Raghunathan}).

Both $\calg_{\Z}=\SL(n,\Z)$ and $\calg_{\Z}=\Sp(n,\Z)$
have the CSP for $n > 2$. With the CSP in play,
an arithmetic subgroup $H$ can be handled using our congruence
homomorphism methodology.
Once the level of a PCS in $H$ is known, most
of the calculations are thereby transferred to finite 
quotients modulo some $\Gamma_{n,m}$. 
So we will be computing with groups 
over integer residue class rings $\Z_m:=\Z/m\Z$
rather than finite fields.

\subsection{Structure of arithmetic subgroups}

\subsubsection{Decidability and principal congruence subgroups}
\label{DecidabilityAndPCS}

Unless stated otherwise, henceforth 
$\calg=\SL(n,\C)$ or $\Sp(n,\C)$ for $n>2$.
Since each arithmetic group in $\calg_{\Q}$ is conjugate
to a group over $\Z$, we confine ourselves to subgroups of $\calg_\Z$
(see \cite[Section~5]{Arithm} for a method to compute a conjugating
matrix; this uses Proposition~\ref{AnotherIntegralityCriterion}).
As we know, if $H$ has finite index in $\calg_\Z$ then 
 $H$ contains some $\Gamma_{n,m}$.
We write $\Gamma_n$ for $\Gamma_{n,1}=\SL(n,\Z)$ or $\Sp(n,\Z)$.

Let $R$ be a commutative unital ring. A \emph{transvection} 
$t\in \mathrm{SL}(n,R)$ is a unipotent
matrix such that $1_n-t$ has rank $1$. 
Denote by $t_{ij}(m)$ the transvection $1_n+e_{ij}(m)$, where 
$e_{ij}(m)$ has $m$ in position $(i,j)$ and zeros elsewhere. Define
\[
E_{n,m} = \langle \hspace{.1mm} t_{ij}(m)\, :  \, i\neq j, \, 1 \leq
i, j \leq n \hspace{.1mm} \rangle
\]
if $\Gamma_n=\SL(n, R)$, and
\begin{align*}
E_{n,m} = & \ \,  \{ t_{i,s+j} (m)\hspace{.75pt}
t_{j,s+i}(m), \hspace{.25pt} t_{s+i,j}(m)\hspace{1pt} t_{s+j,i}(m)
\; | \; 1\leq i < j\leq
s\}\\
& \ \ \cup \hspace{.5pt} \{ t_{i,s+i}(m), t_{s+i,i}(m) \; | \;
1\leq i\leq s\}
\end{align*}
if $\Gamma_n=\Sp(n, R)$ where $n=2s$. The $E_{n,m}$ are
\textit{elementary subgroups of $\, \Gamma_n$ of level $m$}.
Note that $E_{n,1}=\Gamma_n$.

Each arithmetic subgroup in $\Gamma_n$ contains a unique 
maximal PCS. We define the \emph{level} $M(H)$ of an arithmetic 
group $H$ to be the level of its maximal PCS. 
For $\Gamma_n= \SL(n,\Z)$ and $n\geq 3$, the normal closure 
$E_{n,m}^{\hspace{.2pt} \Gamma_n}$ is  
$\Gamma_{n,m}$~\cite[Proposition~1.6]{Arithm}. Similarly,
$E_{n,m}^{\hspace{.2pt} \Gamma_n}$ is the PCS of level
$m$ in $\Gamma_n= \abk \Sp(n,\Z)$ if 
$n> \abk 2$~\cite[Proposition~13.2]{Bass}. 
So at least we have normal generators for a PCS. 
When $\Gamma_n=\SL(n,\Z)$, a full generating set of 
$\Gamma_{n,m}$ is returned by ${\tt GeneratorsPCS}(m)$, 
which encodes the formula in \cite{Venka}. 
The size of this generating set depends only on $n$ (i.e., 
not on $m$). Although a minimal generating set of an arithmetic 
subgroup can be arbitrarily large (see \cite{Venka} again), 
each arithmetic subgroup in $\SL(n,\Z)$ has a $2$-generator 
finite index subgroup~\cite{Meiri}.

Just as we do not need a full generating set of a congruence subgroup to 
compute with a solvable-by-finite linear group, to compute with 
an arithmetic group $H\leq \Gamma_n$ we do not need a full generating 
set for its maximal PCS. 
All we need is the level of $H$.
In the first instance, decidability is then implied by the CSP.
\begin{proposition}[{\cite[Proposition~2.3]{Density}}]
\label{Decidable}
Computing the level of $H$ is decidable.
\end{proposition}

\begin{corollary}[{\cite[Corollary~2.4]{Density}}]
\label{ExplicitlyGiven}
Testing membership of $g\in \Gamma_n$ 
in $H$ and computing an upper bound on $|\Gamma_n:H|$ are decidable.
\end{corollary}

By Corollary~\ref{ExplicitlyGiven}, 
arithmetic subgroups $H\leq 
\Gamma_n$ are explicitly given in the sense of \cite{GSI}.
 Hence the algorithmic problems in
\cite{GSI} for such $H$ are decidable. 

\subsubsection{Computing with subgroups of $\GL(n,\Z_m)$}
\label{ComputingWithGroupsOverZm}

At the moment, algorithms 
for matrix groups over finite rings are less sophisticated than those 
for groups over finite fields. 
Below we sketch the approach in \cite[Section~2]{Arithm} and 
\cite[Section~2]{Density} 
to computing in $\GL(n,\Z_m)$. 

Let $m= p_1^{k_1} \cdots p_t^{k_t}$ where the $p_i$ are distinct
primes and $k_i\geq \abk 1$. By the Chinese Remainder Theorem,
$\chi:\Z_m \rightarrow \Z_{p_1^{k_1}} \oplus \cdots \oplus
\Z_{p_t^{k_t}}$ defined by $\chi(a) = (a_1, \ldots , a_t)$, where
$0\leq\abk a \leq \abk m-1$, $0\leq \abk a_i \leq \abk
p_i^{k_i}-1$, and $a_i\equiv \abk a \mod p_i^{k_i}$, is a ring isomorphism.
\begin{lemma}
\label{FundamIsom}
The map $\chi$ extends to an isomorphism 
\[
\mathrm{Mat}(n,\Z_m)\rightarrow 
{\large \oplus}_{i=1}^t \mathrm{Mat}\big(n,\Z_{p_i^{k_i}}\big)
\] 
which restricts to an isomorphism
\[
\mathrm{GL}(n,\Z_m)
\rightarrow {\textstyle \Pi}_{i=1}^t \mathrm{GL}\big(n,\Z_{p_i^{k_i}}\big).
\]
\end{lemma}

\begin{lemma}\label{PCSFiniteRings}
For $i\geq 1$, let
$K= \{ h \in \GL(n,\Z_{p^k}) \; |\; h \equiv 1_n \mod p^{k-1}\}$. 
\begin{itemize}
\item[{\rm (i)}] 
$K$ is a $p$-group, the PCS of $\GL(n,\Z_{p^k})$ of level 
$p^{k-1}$. 
\item[{\rm (ii)}] $\GL(n,\Z_{p^k})/K\cong \GL(n,p)$.
\end{itemize}
\end{lemma}

Our approach rests on Lemmas~\ref{FundamIsom} and 
\ref{PCSFiniteRings}. First we reduce to $\GL(n,\Z_{p^k})$.
The second part involves computing with
finite $p$-groups and subgroups of $\GL(n,p)$, for which there is
an extensive apparatus~\cite[Sections~7.8 and 9.4]{HoltBook}.
More finely-tuned techniques may be needed;
see, e.g., 
\cite[Section~2]{Density}.

\subsection{Density and computing in arithmetic subgroups}
\subsubsection{Strong approximation}
\label{StrongApproximation}

To compute with arithmetic and dense groups we replace 
finite approximation by
\emph{strong approximation}.

The previous algorithms for solvable-by-finite groups took 
one congruence homomorphism at a time. 
Now we must work with 
congruence images modulo all maximal ideals of $\Z$, i.e., modulo all 
primes $p\in \Z$. This is an option if we have surjection  of 
$H\leq \mathcal{G}_\Z$ onto 
$\varphi_p(\mathcal{G}_{\Z})$ for all but finitely 
many primes $p$; which of course may not happen with
an arbitrary finitely generated subgroup of $\mathcal{G}_\Z$ and
arbitrary $\mathcal G$.
For example, if $m\geq 5$ then reduction modulo $m$ does not surject 
$\GL(n,\Z)$ onto $\GL(n,\Z_m)$.
On the other hand, $\SL(n,\Z)$ surjects onto $\SL(n,\Z_m)$ when 
$m\geq 2$.  Behind these simple observations lies a deep result, the 
\emph{strong approximation theorem} (SAT)~\cite[Window~9]{LubotzkySegal}.
It provides conditions under which $H\leq \mathcal{G}_{\Z}$ surjects
onto $\varphi_p(\mathcal{G}_{\Z})$ for almost all $p$.
A necessary condition for SAT is that $H$ be Zariski dense in $\calg$.
This explains why $\GL(n,\Z)$ does not surject onto $\GL(n,p)$ for
almost all primes $p$: $\GL(n,\Z)$ is not dense in 
$\GL(n,\C)$ because its Zariski closure consists of $g\in \GL(n,\C)$ such that 
$\mathrm{det}(g)^n=1$~\cite[p.~273]{Rapinchuk}.
Neither does density imply SAT; e.g., $\GL(2,\frac{1}{2}\Z)$ is dense in 
$\GL(2,\C)$, but does not surject onto $\GL(2,p)$ modulo any prime 
$p\equiv 1$ modulo $8$~\cite[p.~273]{Rapinchuk}.
However, dense subgroups of $\Gamma_n$ 
satisfy SAT.
\begin{theorem}[{\cite[Window~9]{LubotzkySegal}}]
\label{PreCor30}
If $H \leq \Gamma_n$ is dense in $\mathcal{G}$ then
$\varphi_p(H) = \varphi_p(\Gamma_n)$ 
for almost all primes $p$.
\end{theorem}
Moreover, we have
\begin{theorem}[\cite{Lubotzky97} and {\cite[p.~396]{LubotzkySegal}}]
\label{Cor30}
$H\leq \Gamma_n$ is dense if and only if 
$\varphi_p(H) = \varphi_p(\Gamma_n)$ 
for some prime $p > 3$.
\end{theorem}
Let $\Pi(H)$ be the set of primes $p$ such that 
$\varphi_p(H) \neq \SL(n,p)$.
By Theorem~\ref{PreCor30}, $\Pi(H)$ is finite when $H$ is dense. 

Suppose that $H\leq \Gamma_n$ is arithmetic of level $M$. Then
$\varphi_p(H) = \varphi_p(\Gamma_n)$ for any prime $p$ coprime to $M$
($\varphi_p(\Gamma_n)$ is generated by transvections, and $H$ contains the 
elementary group $E_{n,M}$). Remarkably, the converse (with a tiny number 
of exceptions) is true as well.
\begin{theorem}[{\cite[Section~2]{Density}}]\label{KeyDensity}
\label{PiHEqM}
Let $H\leq \Gamma_n$ be arithmetic of level $M$. 
Then $\varphi_p(H) = \varphi_p(\Gamma_n)$ if and only if 
$p\nmid M$; unless $n=3$ or $4$, $\Gamma_n = \SL(n,\Z)$, 
$\varphi_2(H) = \varphi_2(\Gamma_n)$, and 
$\varphi_4(H)\neq \varphi_4(\Gamma_n)$. In the latter event, 
$M$ is even.
\end{theorem}

Theorem~\ref{KeyDensity} pinpoints the set of prime divisors of the level 
of arithmetic $H\leq \Gamma_n$: barring the exceptions in
Theorem~\ref{PiHEqM}, it is exactly $\Pi(H)$. This set
is input for our algorithm to compute $M$;
see Subsection~\ref{ComputingMEtc} and \cite[Section~2.4]{Density}.

We can compute with dense rather than merely arithmetic groups 
in $\Gamma_n$. Each dense subgroup $H$ has a unique minimal arithmetic 
overgroup $\mathrm{cl}(H)$the intersection of all arithmetic groups 
in $\Gamma_n$ that contain $H$. The \emph{level} of $H$ is then defined 
to be the level of $\mathrm{cl}(H)$. 
We implemented an algorithm to compute 
this  `arithmetic closure' $\mathrm{cl}(H)$~\cite[Section~3.3]{Density};
see Subsection~\ref{ArithmeticAlgorithms}. 
Hence, we can investigate $H$ by applying 
algorithms for arithmetic groups to $\mathrm{cl}(H)$. Perhaps 
$\mathrm{cl}(H) = \Gamma_n$. This does happen: 
see \cite{SHumphries} for examples of free subgroups
$H$ of $\SL(n,\Z)$, $n>2$, generated by $n$ transvections, that surject 
onto $\SL(n,p)$ modulo all primes $p$.
Since $\SL(n,\Z)$ is virtually free only for $n=2$, $H$ has infinite index
in $\SL(n,\Z)$, i.e., is thin. 
Note that, for $n\geq 5$, the only arithmetic subgroup in $\Gamma_n$ 
that surjects onto $\varphi_p(\Gamma_n)$ for all primes $p$   
is $\Gamma_n$ itself~\cite[Corollary~2.14]{Density}.

\subsubsection{Density testing and SAT}
\label{DensityTestingSubsection}
While decidability of arithmeticity testing is unknown, 
we can test density, and this
serves as an initial check along the way to 
settling arithmeticity of individual examples. 

A Monte Carlo algorithm to test density is given in \cite{RivinIII}. 
It uses
\begin{theorem}
Suppose that $H\leq \Gamma_n$
contains non-commuting elements $g_1$, $g_2$ such that the Galois 
group of the characteristic polynomial of $g_1$ is $\mathrm{Sym}(n)$, 
and $g_2$ has infinite order. Then either $H$ is dense, or
$\calg=\Sp(n,\C)$ and the closure of $H$ over $\C$ is the product of 
$n/2$ copies of $\SL(2,\C)$.
\end{theorem}

There are intrinsic {\sf GAP} procedures to compute Galois groups and 
test equality with $\mathrm{Sym}(n)$. 
And finiteness can be tested over any field 
(Section~\ref{FinitenessSection}).
Elements of $\Gamma_n$ 
with associated Galois group equal to the symmetric group are 
\emph{generic}: a `random' element of $\Gamma_n$ is likely to 
satisfy the criteria~\cite[Theorem~1.4]{RivinIII}. Algorithm~1 
of \cite{RivinIII} is 
Monte Carlo, and may incorrectly report that input is not dense; 
but the probability of error is small due to the abundance of 
elements $g_1$ and $g_2$.

Another density testing algorithm in \cite{RivinIII} accepts a finitely 
generated subgroup $H$ of semisimple $\calg$ 
in characteristic zero. It uses the fact that $H$ is dense if 
and only if (i)~$H$ is infinite, and (ii)~the adjoint representation 
$\mathrm{ad}(H)$ on the Lie algebra of $\calg$ is absolutely
irreducible. If $\calg= \abk \SL(n,\C)$ then the Lie algebra 
consists of all matrices with zero trace, so has dimension $n^2-1$. 
If $\calg= \abk \Sp(n,\C)$ then the algebra has dimension $(n^2+n)/2$, and 
consists of all matrices of the form
${\tiny 
\left(\renewcommand{\arraycolsep}{.12cm} 
\begin{array}{cc} A & B \\
C & A^\top
\end{array} \! \right)}$
where $B$ and $C$ are symmetric.
We construct $\mathrm{ad}(H)$ from a finite generating set for
$H$. Checking absolute irreducibility is routine; cf.~\cite[p.~404]{Tits}. 
So we have deterministic density testing, too~\cite[Section~5]{Density}.

We get a one-way density test from Theorem~\ref{Cor30}, i.e., if 
$\varphi_p(H) = \varphi_p(\Gamma_n)$ for some prime $p>3$ then $H$ is
dense. This raises the problem of realizing SAT computationally: 
for dense $H$, compute the set $\Pi(H)$.
The problem is addressed in \cite[Section~3.2]{Density}
and in \cite{DensityFurther}. We label the procedures from
\cite{Density,DensityFurther} for computing $\Pi(H)$ as 
${\tt PrimesForDense}(H)$. 

\subsection{Algorithms for computing with dense and arithmetic groups}
\label{ArithmeticAlgorithms}

In this subsection we outline procedures for dense 
(including arithmetic) groups in $\Gamma_n$, $n>2$.

\subsubsection{Computing the level and related procedures}
\label{ComputingMEtc}

We tailor congruence homomorphism methods to 
properties of arithmetic or dense groups (CSP, SAT). The core  
part is computing the level.

$\tt LevelMaxPCS$ accepts a finite set $S\subseteq \Gamma_n$ that 
generates a dense group $H$, and returns the level of 
$H$. To give a flavor of the computation, we 
quote a technical result (`stabilization lemma') from \cite{Density}. 
For $n>2$ and $H\le \Gamma_n$, set
\begin{equation*}
\delta_H(m)=|\Gamma_n:\Gamma_{n,m}H|; \tag{$\diamond$}
\end{equation*}
i.e., $\delta_H(m)=|\varphi_m(\Gamma_n):\varphi_m(H)|$.
\begin{lemma}[{\cite[Lemma~2.16]{Density}}]
\label{deltaGen} \ 
\begin{itemize}
\item[{\rm (i)}]
Suppose that $\delta_H(kp^{a})=\delta_H(kp^{a+1})$ for some prime 
$p$, positive integer $a$, and $k$ coprime to $p$.
Then $\delta_H(kp^{b})=\delta_H(kp^{a})$ for all $b\geq a$.
\item[{\rm (ii)}]  Let $p$, $a$, and $k$ be as in {\rm (i)}. 
Then $\delta_H(lp^b)=\delta_H(lp^{a})$ for all $b\geq a$ and any
multiple $l$ of $k$ such that $\pi(l)=\pi(k)$.
\end{itemize}
\end{lemma}

$\tt LevelMaxPCS$ is recursive, and embodies theory of dense groups 
as in Subsection~\ref{StrongApproximation}.
Its input includes ${\tt PrimesForDense}(H)$. This gives the 
primes dividing the level of $H$ (with some exceptions as in 
Theorem~\ref{PiHEqM}, which cause no difficulty in the computation).
The highest power of each prime $p$ in the prime factorization of the 
level is determined. Since by $(\diamond)$ all $\delta$-values 
encountered in its recursion loop 
are bounded, the procedure terminates. 

$\tt LevelMaxPCS$ solves the following problems.
We can compute $|\Gamma_n:H|$ if $H$ is arithmetic
(cf.~Lemma~\ref{deltaGen}). Membership testing is easy: 
$g\in \Gamma_n$ is in $H$ if and only if $\varphi_M(g)\in \varphi_M(H)$,
where $M= \abk {\tt LevelMaxPCS}(S)$.
This extends to ${\tt IsSubgroup}(H,L)$, which tests whether 
$L\leq \Gamma_n$ is contained in $H$. When $H$ is arithmetic, 
$\tt LevelMaxPCS$ returns $1$ if and only if  $H=\Gamma_n$.

\subsubsection{Subnormality}

The contrast between linear groups over rings and groups over 
fields is stark when we look at their subnormal subgroups.
We touched on this point in Subsection~\ref{CSPSection}  
and elaborate on it now.

In this subsection, $\Gamma_n = \SL(n,\Z)$.
Statements can be duly modified for $\Gamma_n = \Sp(n,\Z)$.
As in \cite[p.~166]{JW}, for $h=(h_{ij})\in \Gamma_n$ we denote by $\ell(h)$ 
the ideal of $\Z$ generated by
\[
\{ h_{ij} \st i\neq j, \, 1\leq i , j \leq n\} \cup \{
h_{ii}-h_{jj} \st  1\leq i , j \leq n \},
\]
i.e., the lcm of the non-diagonal entries of $h$ and differences of 
the diagonal entries. 
Then $\ell(A) :=\abk \sum_{a\in A} \ell(a)$ for $A\subseteq
\abk \Gamma_n$. Let $Z_{n,m}$ denote the full preimage of the center (scalar
subgroup) of $\GL(n,\Z_m)$ in $\Gamma_n$ under $\varphi_m$. 
So $\ell(A)$ is the smallest ideal $m\Z$ such that
$A\subseteq \abk Z_{n,m}$. We define $\ell(A)$ unambiguously as the non-negative
integer modulo $m$ that generates $\ell(A)$; e.g., $\ell(Z_{n,k})
= \abk \ell(\Gamma_{n,k}) = \abk k$.

If $H=\gpess\leq \GL(n,\Z)$ then 
$\ell(H) = \ell(S)$~\cite[Lemma~1.22]{Arithm}.
This gives a procedure ${\tt El}(S)$ that returns $\ell(H)$ for
$H=\langle S\rangle$. 

We need the following `sandwich theorem'.
\begin{theorem}[\cite{JW}]
$H\leq \GL(n,\Z)$ is subnormal if and only if
$\Gamma_{n,k^e} \leq H \leq Z_{n,k}$
for some $k$, $e$.   If $|\GL(n,\Z):H|<\infty$ then $\Gamma_{n,M}$ is the 
maximal PCS of $H$, in which case $H$ is subnormal in $\GL(n,\Z)$ if and 
only if $M$ divides $\ell(H)^t$ for some $t$.
\end{theorem}
\begin{corollary}[{\cite[Corollary~1.28]{Arithm}}]
$H\unlhd \GL(n,\Z)$ is subnormal if and only if $\ell(H)$ is the level of $H$.
\end{corollary}
So the normal closure 
of an arbitrary subgroup $H$ of $\GL(n,\Z)$ 
is $\langle H,\Gamma_{n,\ell(H)}\rangle$.
The procedure ${\tt NormalClosure}$ accepts (finite) $S$ 
for $H=\langle S\rangle$, computes $\ell = {\tt El}(S)$, and returns
the union of $S$ and ${\tt GeneratorsPCS}(\ell)$.

For arithmetic $H$, ${\tt IsSubnormal}$
returns $\tt true$ if there is non-negative $e\in \Z$ such that 
$M= {\tt LevelMaxPCS}(S)$ divides ${\tt El}(S)^e$. 
If ${\tt El}(S) =1$, i.e., $M={\tt El}(S)$,
then the procedure ${\tt IsNormal}(S)$ tests whether $H$ is normal
in $\GL(n,\Z)$. We compute the normalizer $N_{\Gamma_n}(H)$ of arithmetic 
$H\leq \Gamma_n$ as the full preimage of the normalizer of a 
congruence image.

\subsubsection{The orbit-stabilizer problem}
\label{OSProblem}

Let $G\leq \glnf$, and $u$, $v\in \F^n$.
The \emph{orbit-stabilizer problem} is:

\vspace{7pt}

 (OP) \, Decide whether there 
is $g \in G$ such that 
$gu=v$;
if so, find $g$.

\vspace{2.5pt}

(SP) \,  Compute a  
generating set of $\mathrm{Stab}_G(u)
\! = \! \{ g \in G  \; |\;  gu = u\}$.

\vspace{7pt}

\noindent 
The orbit problem (OP) is related to two other fundamental
algorithmic questions,
membership and conjugacy testing (see \cite[3., p.~239]{Dixon85}).
Since the conjugacy problem can be undecidable in $\SL(4,\Z)$, 
(OP) can be undecidable. And we only attempt to solve the stabilizer 
problem (SP) knowing beforehand that $\mathrm{Stab}_G(u)$ is finitely 
generated 
(see \cite[6., p.~239]{Dixon85} for a simple example where it is not).

The above difficulties steer us toward arithmetic groups 
$H$ in $\Gamma_n= \SL(n,\Z)$. 
The orbit-stabilizer problem is decidable for any 
explicitly given arithmetic subgroup of $\calg$~\cite{GSI};
this implies that it is decidable for $H$. A practical algorithm 
is proposed in \cite[Section~4]{Arithm}. 
We first solve the orbit-stabilizer problem 
for subgroups of $\GL(n,\Z_m)$ acting on $\Z_m^n$, then solve 
it for $\Gamma_{n,m}$ acting on $\Z^n$, and patch together 
the two solutions. The  second stage of the method 
uses the next theorem (see
\cite[Proposition~4.10 and Theorem~4.13]{Arithm}). 
\begin{theorem}
Non-zero elements $u = (u_1,\ldots , u_n)^\top$ and 
$v = (v_1,\ldots , v_n)^\top$  of $\Z^n$ are in the same 
$\Gamma_{n,m}$-orbit if and only if 
\begin{itemize}
\item the $u_i$ generate
the same ideal $a\Z$ of $\Z$ as the $v_i$, and 
\item $u_i \equiv v_i \, \mathrm{mod} \, am$ for
$1\leq i\leq n$. 
\end{itemize}

\vspace{.5pt}

\noindent In particular, if $m=1$ then $u$, $v$ are in the same 
$\Gamma_n$-orbit if and only if the $u_i$ generate
the same ideal  as the $v_i$.
\end{theorem}
A similar result is true over $\Z_m$, for use in the first stage.
\begin{proposition}[{\cite[Proposition~4.7]{Arithm}}]
Non-zero elements $u = (u_1,\ldots , u_n)^\top$ and 
$v = (v_1,\ldots , v_n)^\top$ of $\Z_m^n$ are in the same 
$\SL(n,\Z_m)$-orbit if and only if the $u_i$ generate
the same ideal of $\Z_m$ as the $v_i$.
\end{proposition}

So we obtain a procedure ${\tt Orbit}(S,u)$ that returns a solution of (OP)
for $H=\langle S\rangle\leq \SL(n,\Z)$ and $u\in \Z^n$, $n>2$.

By \cite[p.~744]{GSI}, $\mathrm{Stab}_{\Gamma_n}(u)$ is finitely generated:
it is conjugate to the affine group of matrices of the form 
\[
 {\small \left(\renewcommand{\arraycolsep}{.15cm} 
\begin{array}{cc} 1 & *
\\
\vspace*{-9pt}& \\
0 & \Gamma_{n-1}
\end{array} \! \right)}.
\]
Let $\sigma = {\tt Orbit}(\Gamma_n,u)$. Remember
that we know a generating set of $\Gamma_n$.
Then $\mathrm{Stab}_{\Gamma_n}(u)$ is generated by the 
$\sigma$-conjugates of 
$t_{12}(1),\dots , t_{1n}(1)$ and $\mathrm{diag}(1,x)$ as $x$ 
runs over a generating
set of $\Gamma_{n-1}$.
Similarly, $\mathrm{Stab}_{\Gamma_{n,m}}(u)$ is generated by the 
$\sigma$-conjugates of $t_{12}(m),\dots ,t_{1n}(m),\mathrm{diag}(1,x)$ 
as $x$ runs over a generating set of $\Gamma_{n-1,m}$ (as in 
\cite{Venka}). This completes solution of (SP) for 
$G\leq \SL(n,\Z_m)$ acting on $\Z_m^n$, 
hence solving the orbit-stabilizer problem for finitely 
generated subgroups of finite index in $\Gamma_n$.

\subsection{Experiments}\label{ExperimentalNonSF}

The algorithms from this section are joint work with Alexander Hulpke, 
and have been implemented in {\sf GAP}. Below is a 
small sample of experiments conducted with the software (see \cite{GAPDoc}).

One set of test groups come from low-dimensional topology.
These were evidence supporting the conjecture of 
Lubotzky on $2$-generator arithmetic subgroups, verified 
in \cite{Meiri} (see Subsection~\ref{DecidabilityAndPCS}). Let 
\[
\mathcal{C} = \langle x, y, z\; | \; 
z  x z^{-1} = x y, \, z y z^{-1} = y x
 y \rangle, 
\]
the fundamental group of the figure-eight knot complement. 
A representation $\beta_T:\mathcal C\rightarrow \SL(3,\Z)$ for 
$T\in \Z\setminus \{0\}$ is given in \cite{LongReidI}, where 
$H_T=\langle \beta_T(x), \beta_T(y) \rangle$ is shown to be
arithmetic. Finding $|\SL(3,\Z):H_T|$ was an open problem 
in \cite{LongReidI}.   
We can compute the level and the index using 
${\tt PrimesForDense}$ and 
$\tt LevelMaxPCS$; see \cite[Section~6]{Arithm} and 
\cite[Section~4]{Density}.
For example, it took 892.6 seconds to find 
the level $2^{7}5^{6}29 {\cdot}67 {\cdot}193$ and index 
$2^{42}3^{5}5^{25}7^{4}13{\cdot}31^{2}67{\cdot}1783$
for $T= 100$.

Another batch of test groups comes from an intermingling of group theory,
differential equations, and theoretical physics. Let
$G(d,k) = \langle U, T\rangle$ where
\[
U =  
{\footnotesize \left(\begin{array}{crrr} 1 & 1 & 0 & \ 0
\\
0 & 1 & 0 & \ 0
\\
d & d & 1 & \ 0
\\
0 & -k & -1 & \ 1
\end{array} \right)},
\quad  T =  
{\footnotesize
\left(\begin{array}{cccc} 
1 & 0 & 0 & 0
\\
0 & 1 & 0 & 1
\\
0 & 0 & 1 & 0
\\
0 & 0 & 0 & 1
\end{array} \right)}.
\]
There are $14$ pairs $(d,k)$ such that 
$G(d,k) \leq \mathrm{Sp}(4,\Z)$ is the monodromy group of a
 hypergeometric ordinary differential equation
associated to Calabi-Yau threefolds. Of those, seven 
are arithmetic and the rest are thin~\cite{SinghVenky}.
We can study $G(d,k)$ by first locating an arithmetic group
in $\Sp(4, \Z)$ containing $G(d, k)$; 
see \cite{MonodromyAppendix}. Our procedures 
complete all computations quickly, and find
the minimal arithmetic overgroup of $G(d,k)$
(some of these for the first 
time)~\cite[Section~4.2]{Density}.

\section{Where to next?}

We have used our methodology for computing with finitely generated 
linear groups to solve a variety of algorithmic 
problems. Other problems await a similarly satisfactory treatment. 

\subsection{Solvable-by-finite groups}

One of the foremost open computational problems for 
solvable-by-finite groups $G\leq \GL(n,\F)$ 
is subgroup membership testing. 
This problem is decidable for $\F= \Q$~\cite{Kopytov2},
so it is decidable for finitely generated linear groups of finite 
rank (Proposition~\ref{SFSufficesForQLinear}).
However, as yet we do not have a practical algorithm.
One hurdle to overcome is that the unipotent radical $U(G)$ may not 
be finitely generated.
However, $U(G)$ has finite rank (and is polyrational), so we could 
attempt to design an algorithm using methods 
from Subsections~\ref{InStructureSF} and \ref{NewFRSection}. 
This would amount to a 
computational realization of Hall's theory of infinite nilpotent 
groups, including calculation of Hall polynomials.
With an eye on practicality, we note another technique: 
replace computing in the torsion-free nilpotent group $U(G)$ by 
computing in related Lie algebras (see, e.g., \cite{AssmannMalcev}).

Partial classes of solvable-by-finite linear groups offer fresh 
opportunities. Constructing (faithful) linear representations of 
polycyclic-by-finite groups falls under this heading.
An algorithm to solve this problem for finitely generated torsion-free 
nilpotent groups is given in \cite{Nickel}. More ambitiously, we 
would seek an algorithm to construct a representation over $\Q$  
of a finitely generated finite rank subgroup of $\glnf$. 

\subsection{Non-solvable-by-finite groups}\label{NonSFAndFurtherOP}

Let $G$ be a finitely generated non-solvable-by-finite subgroup of $\glnf$.
The procedure $\tt IsSolvableByFinite$ from 
Subsection~\ref{TAltRealization} only attests to the existence 
of a free group in $G$, without constructing one. Taking
advantage of the ubiquity of free subgroups in non-solvable-by-finite 
linear groups~\cite{Aoun,Epstein},
we could try picking a `random'  finitely generated subgroup and testing 
whether it is free (here `random' has several interpretations~\cite{RivinPick}). 
This runs into the open problem of testing whether finitely generated 
linear groups are free. 
The difficulty of freeness testing is evinced by 

\vspace{-4pt}

\[
H(x) := \Big\langle 
{\tiny \left( \begin{array}{cc}
1 & x \\
0 & 1\end{array}
\right), 
\left( \begin{array}{cc}
1 & 0 \\
x & 1\end{array}
\right)} \Big\rangle \leq  \GL(2, \mathbb{C}).
\]

\vspace{2.5pt}

\noindent
This group is free for $|x| \geq 2$, 
and for some $x$ such that $|x| < 2$.
However, often it is still not known whether $H(x)$ is free
(see \cite{Beardon}). 

Construction of `large' free subgroups would be 
interesting. For a non-solvable-by-finite subgroup $H$ of an
algebraic group $\mathcal G$ (say 
$\mathcal{G} = \overline{H}$), `large' means Zariski dense. If $n > 2$,
$\mathcal G = \Sp(n,\C)$ or 
$\SL(n,\C)$, and $H = \mathcal G_\mathbb{Z}$,  
then such free subgroups are plentiful
and provide examples of thin matrix groups. 
A related problem is construction of strongly dense subgroups of 
$\mathcal G$, i.e., free dense subgroups $H$ in which each non-cyclic 
subgroup is again dense. 
See \cite{BrGGT} for insights on the importance of strongly dense 
subgroups, their existence, 
and a link to the Banach--Tarski paradox.

The group $\SL(3, \Z)$ is a source of exciting open problems.
As usual, subgroup membership testing is crucial: it is not known whether 
this problem is decidable in $\SL(3, \Z)$.
We also mention (i)~the Howson property: 
is the intersection of two finitely generated subgroups of $\SL(3, \Z)$ 
finitely generated? (ii)~coherence: is every finitely generated subgroup 
of $\SL(3, \Z)$ finitely presentable? Computer experimentation 
may help to answer these questions; cf.~\cite{LongReidI} and 
\cite[Section~4.1]{Density}. The questions have a negative answer 
in degrees greater than $3$ and a positive answer in degree $2$. 

We pointed to arithmeticity testing as a significant open problem; 
cf.~Section~\ref{522} and \cite{Sarnak}.   
Constructing presentations of (finitely presentable) arithmetic groups  
can be hard; a practical algorithm even for finite index 
subgroups of $\SL(n,\Z)$ is lacking. 
Recent progress is \cite{NebeRecent}, which gives algorithms to compute 
a presentation of the unit group of an order in a semisimple $\Q$-algebra.
The problem of computing linear 
representations of finitely generated abstract groups
resurfaces (especially 
finitely presented groups that are not necessarily 
solvable-by-finite~\cite{PlesSur}).
Here we  note the achievements of \cite{LongReidII,LongLatest}
in constructing faithful representations of triangle groups.

\subsection*{Acknowledgments}
A.~S.~Detinko was supported by a
Marie Sk\l odowska-Curie Individual
Fellowship grant H2020 MSCA-IF-2015, No.~704910 (EU Framework Programme 
for Research and Innovation). 
D.~L.~Flannery was funded by the Irish Research Council 
(New Foundations 2015).

\bibliographystyle{amsplain}
\def\cprime{$'$}
\providecommand{\bysame}{\leavevmode\hbox to3em{\hrulefill}\thinspace}
\providecommand{\MR}{\relax\ifhmode\unskip\space\fi MR }
\providecommand{\MRhref}[2]{%
  \href{http://www.ams.org/mathscinet-getitem?mr=#1}{#2}
}
\providecommand{\href}[2]{#2}

\bigskip

\end{document}